%% file: Balrog_Wrapper.tex
%
%

\def\wrapperUsed{WrapperPresent}
\magnification=\magstep1

\footline={\relax}

\
\vskip2.4truein

\centerline{\bf Inconsistency of Primitive Recursive Arithmetic}
\bigskip
\centerline{\bf Edward Nelson}

\bigskip
\bigskip
\centerline{With an}
\bigskip
\centerline{Introduction by Sarah Jones Nelson}
\bigskip
\centerline{ and an }
\bigskip
\centerline{ Afterword by Sam Buss and Terence Tao }
\vfill
\vfill
\eject

\input SarahIntro_arXiv.tex

\vfill\eject

\input NelsonPhotoPage.tex

\pageno=1

\footline={\hfill \the\pageno\hfill}

\input Balrog_Main_arXiv.tex

\vfill\eject

\footline={\relax}\headline={\relax}

\input BussTao_Afterword_arXiv.tex

\vfill \eject
\input Appendix13i.tex

\end

%% file: SarahIntro_arXiv.tex
\font\sambigrm=cmr10 scaled \magstep2

{
\ifx\wrapperPresent\relax
\magnification=\magstep1
\fi
\hsize=5.0in
\hoffset=0.25in

\
\vskip1cm

\centerline{
\sambigrm  Introduction}

\vskip1cm

In September of 2011 Edward Nelson announced that he had a proof of the inconsistency of Peano Arithmetic. He had devoted twenty-five years to constructing the proof. When Terence Tao and Daniel Tausk independently found an error, Ed withdrew his claim at once and cheerfully returned to work on it the next day. By March of 2013 -- confident that he had corrected the error -- he wrote a project proposal with ``a crucial new insight: technically, bounds in the Hilbert-Ackermann consistency theorem depending only on rank and level, but not on the length of the proof.'' The project he proposed ``challenges the entire current understanding and practice of mathematics....It will radically change the way mathematics is done. This will affect the philosophy of mathematics, how the nature of mathematics is conceived. It will also affect the sciences that use mathematics, especially physics....It addresses a very big question: does mathematics consist in the discovery of truths about some uncreated eternal reality (the traditional Platonic view), or is it a humble human endeavor to construct abstract patterns that will be sound, free of all contradiction, in the hope that some will be beautiful, uplifting the human spirit, and that some (not necessarily different ones) will be of practical use to improve the human lot?''

The following excerpts from Ed's proposal describe his vision of a new mathematics and the open question of the consistency of Peano Arithmetic.
\vskip 1em

``Peano Arithmetic is one of the simplest and most fundamental of mathematical
theories. Its consistency, however, has not been proved by any means that all
mathematicians accept. It implies that all primitive recursive functions are
total, but they are directly defined only for numerals, and the argument that
the values always reduce to numerals is circular. The proposal is to complete a
proof that Peano Arithmetic is in fact inconsistent. The principal output will
be a book entitled `Elements'. The outcome will be a major change in the way
mathematics is done, with philosophical and scientific consequences.

``The guiding spirit of this investigation is that mathematics is not some
uncreated abstract reality that we can take for granted and explore, but that
it is a human endeavor in which one should begin by looking at the very
simplest concepts without taking them for granted.

``Numbers are constructed from $0$ by successively taking successors; ${\rm S
\ldots S}0$ is called a numeral. Definitions of primitive recursive functions,
such as addition, multiplication, exponentiation, superexponentiation, and so
forth, are schemata for constructing numerals. They define a value for $0$ and
then a value for ${\rm S}x$ in terms of the value for $x$. But when numerals
are substituted for the variables in such a schema it is not clear that it
defines a numeral: the putative number of steps needed to apply the definitions
can only be expressed in terms of the expressions themselves. The argument is a
vicious circle. Consequently, the consistency of Primitive Recursive
Arithmetic, and {\it a fortiori} of Peano Arithmetic ($\rm P$), is an open
question.

``Here is a nontechnical description of how I propose to show that $\rm P$ is
inconsistent. We start with a weaker theory $\rm Q$ and by relativization
techniques extend it to a stronger theory $\rm Q^*$. Proofs in $\rm Q^*$ reduce
to proofs in $\rm Q$. $\rm Q^*$ arithmetizes $\rm Q$ itself -- that is, it
expresses the syntax of $\rm Q$ by a term of $\rm Q^*$ that here I shall denote
by $\rm @Q$. Remarkably, $\rm Q^*$ proves that there is no open proof of a
contradiction in $\rm @Q$. (`Open' means that the proof has no quantifiers,
i.e., symbols for `there exists' and `for all'.) All of this was done in my
book `Predicative Arithmetic' (Princeton University Press, 1986) and is being
redone with complete proofs in the book `Elements', a work in progress that is
the subject of this proposal. The Hilbert-Ackermann consistency theorem implies
that there is no proof of a contradiction, even with quantifiers, in $\rm @Q$.
This theorem can be only partially established in $\rm Q^*$. The crucial new
insight is that it can be proved provided there are bounds on features of the
proof (called rank and level) but emphatically not depending on the length of
the proof. These bounds on rank and level cannot be proved in $\rm Q^*$, but for
each specific proof, $\rm P$ proves that $\rm Q^*$ proves them!  This opens the
way to exploit the stunning proof without self-reference of G\"{o}del's second
incompleteness theorem by Kritchman and Raz (Notices of the American
Mathematical Society, December 2011). The upshot is that $\rm P$ proves that
$\rm Q$ is inconsistent. But, as is well known, $\rm P$ also proves that $\rm
Q$ is consistent. Therefore Peano Arithmetic $\rm P$ is inconsistent.''
\vskip1em

Ed was the only living mathematician who could argue from purely syntactic reasoning without the traditional semantics established by Plato. John Conway suggested to me that this might explain why no one has fully understood Ed's deeply unique insights into the foundations of contemporary mathematics. He stood alone in the world, courageous as a formalist of a new ontology of integers: proof that completed infinities do not exist and that human minds invented numbers never discovered or revealed from platonic forms of any fundamental reality. I am hopeful that the mathematical community will boldly investigate ``Elements'' and the unshakeable foundations Ed sought to build.

\vskip1cm

\hfill Sarah Jones Nelson
\vskip 1ex

\hfill {\it September 10, 2015}

\vfill

}

\ifx\wrapperUsed\nil
\def\thisEnd{\end}
\else
\def\thisEnd{\relax}
\fi

\thisEnd

%% file: NelsonPhotoPage.tex
\ \vfill

\input epsf
\centerline{\epsfxsize=2in \epsfbox{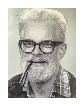}}
\vskip-6pt
\centerline{\hbox to 2in{\hfill \font\sixrm=cmr6 \sixrm Photo courtesy Sarah Jones Nelson}}
\vskip1em
\centerline{\font\twelverm=cmr12 \twelverm Edward Nelson}
\vfill

\eject

%% file: Balrog_Main_arXiv.tex
\magnification=\magstep1
\input NelsonarXivBalrogMacros.tex

\input hyperbasics.tex
\input colordvi.tex

\def\pdfklink#1#2{\hbox{\textBlue{#1}\textBlack}}
\def\pdfKlink#1#2{\hbox{\textBlue{#1}\textBlack}}

\def\smallskip{\medskip}

\centerline{\bf Inconsistency of Primitive Recursive Arithmetic}

\medskip

\centerline{\bf Edward Nelson}

\bigskip

\centerline{\bf //   1. Against finitism}

\medskip

Primitive Recursive Arithmetic (PRA) was invented by Skolem [Sk] in response to \ii Principia Mathematica.\/
with the express purpose of avoiding quantification over infinite domains. His tools were primitive
recursion and induction.

Familiar examples of primitive recursions are
$$\eqalign{
x+0&=x\cr
x+\ro Sy&=\ro S(x+y)\cr
\noalign{\vskip3pt}
x\cdot0&=0\cr
x\cdot\ro Sy&=x+(x\cdot y)\cr
\noalign{\vskip3pt}
x\up0&=\ro S0\cr
x\up\ro Sy&=x\cdot(x\up y)\cr
\noalign{\vskip3pt}
x\Uparrow0&=\ro S0\cr
x\Uparrow\ro Sy&=x\up(x\Uparrow y)\cr
}$$

A numeral is a term containing only S and 0, and a primitive recursive (PR) number is a variable-free term all of
whose function symbols are PR. The finitary credo is that PR~numbers
reduce to numerals by applying the equations a sufficient number of times. If indeed that were so, the
applications used could be counted by a numeral. But in general the number of applications
needed can be expressed only in terms of PR~numbers themselves---the argument is blatantly circular.

The objection being raised here is not some vague semantic ``ultrafinitistic'' assertion that some PR~numbers
are so big they don't really exist. Certainly the PR number
$\rm SS0\Uparrow SSSSSS0$ exists: here it is, in front of our eyes, consisting of eleven symbols.
The problem is syntactical. Let A~be a formula that holds for~0 and is such that whenever it holds for~$x$
it holds for~$\ro Sx$. Then A~holds for any numeral~n;
this follows from the basis $\ro A_x(0)$ by modus ponens applied as many times as there are occurrences of~S
in~n. But the postulation of induction, implying that A~holds for every PR~number, is
an expression of the finitary credo.

PRA is a boldly speculative attempt to treat PR~numbers as if they were equal to numerals.
We shall see that it is inconsistent.

\bigbreak

\centerline{\bf //   2. Outline}

\medskip

The next section describes the notational and terminological conventions used in this paper
(including the present outline) and formulates PRA as a formal system. Section~$//   4$ explicitly defines
a binary function symbol Eq such that $\vdash\'Eq'(x,y)=0\iffF x=y$ and
\S$//   5$ introduces bounded quantifiers. Section~$//   6$ introduces strings and their
combinatorics, and \S$//   7$ uses them to formulate arithmetization.
Section~$//   8$ establishes a form of Chaitin's theorem~[Ch]. The heart of the paper is \S$//   9$,
which constructs a subsystem of PRA that proves the consistency of its own arithmetization.
The final \S$//   10$ exploits the method of the Kritchman-Raz proof
\yy\pdfklink{www.ams.org/notices/201011/rtx101101454p.pdf}{www.ams.org/notices/201011/rtx101101454p.pdf}\zz,
based on Chaitin's theorem and the surprise examination paradox,
together with the self-consistency result of \S$/   9$, to derive a contradiction in PRA.

Numbered and italicized steps occur throughout the paper to serve as a fuller outline.

\bigskip

\centerline{\bf /   3. The formal system PRA}

\medskip

{\it \#1. Formulate the axioms and rules of inference of\/ {\rm PRA}.}

\smallskip

The \ii symbols.\/ of PRA are \ii variables., \ii function symbols.,
$=$, $\Neg$, and~$\vee$. A ``decorated letter'' is a letter with zero or more
digits as subscript and zero or more primes as superscript. We use decorated~s to stand for symbols.
Decorated italic letters are variables, and decorated x~y~z~w stand for variables.
An \ii expression.\/ is a concatenation of symbols; decorated u~v stand for expressions, f~g~h for function
symbols. Each symbol~s has an \ii index.\/ (or arity), denoted by $\rm\iota s$,
specifying how many arguments it takes. A symbol is 0-ary, unary, binary, and so forth, according as its index
is 0, 1, 2, and so forth. A \ii constant.\/ is a 0-ary function symbol; decorated e stands for constants.
Variables are 0-ary, $\Neg$ is unary, = and $\vee$ are binary.
\ii Terms.\/ are defined recursively as follows: x is a term; if $\rm u_1, \ldots, u_{\iota f}$ are terms
then $\rm fu_1\ldots u_{\iota f}$ is a term. Decorated a~b~c~d stand for terms.
An \ii equation.\/ is $\rm {=}ab$. \ii Formulas.\/
are defined recursively as follows: equations are formulas; if u and v are formulas,
so are $\rm\Neg u$ and $\rm\vee uv$. (There are no quantifiers in PRA.) Decorated A~B~C~D~H stand for formulas.
We frequently use infix notation for binary symbols; thus $\rm a=b$ abbreviates $\rm{=}ab$ and $\rm A\vee B$
abbreviates $\rm{\vee}AB$. The use of infix notation often requires groupers to avoid ambiguity; we use
parentheses to group terms, and brackets and braces to group formulas. Some other useful abbreviations
are \hbox{$\rm A\andD B$} for $\rm\Neg[\Neg A\orR\Neg B]$, $\rm A\impP B$ for \hbox{$\rm\Neg A\orR B$}, and
\hbox{$\rm A\iffF B$}
for \hbox{$\rm[\Neg A\orR B]\andD[A\orR\Neg B]$}. The symbol~$\Neg$ binds tightly, and $\orR$
and $\andD$ bind more tightly than $\impP$ and $\iffF$.
Apart from these precedence rules, infix symbols are associated from right to left.
Function symbols are \ii nonlogical symbols..

If $\ell$ is a decorated letter, $\vec{\hskip1pt\ell}$ abbreviates
$\ell_1\ldots\ell_\mu$
for some~$\mu$ called the \ii multiplicity.\/ of~$\vec{\hskip1pt\ell}$ and denoted by $\mu\vec{\hskip1pt\ell}$.
The notation $\rm a\ne b$ abbreviates $\rm\neg a=b$. Let
$\rm u_{\vec x}(\vec a)$, where $\rm\mu\vec x=\mu\vec a$, be the expression, called an \ii instance.\/ of~u,
obtained by replacing each occurrence of~$\rm x_\nu$ in~u by~$\rm a_\nu$, for all $\nu$ with
$\rm1\le\nu\le\mu\vec x$. Whenever we write $\rm u_{\vec x}(\vec a)$ it is understood that
$\rm\mu\vec x=\mu\vec a$.

Sometimes parentheses and commas are inserted into terms to enhance readability.
Although $\rm A_x(x)$ and A are the same, the redundant notation $\rm A_x(x)$ often increases readability.

A \ii truth valuation.\/ on A is a function $\tau$ from the equations in A to $\{\true,\false\}$.
We extend $\tau$ to all subformulas of A, keeping the notation $\tau$, by letting $\rm\tau(\Neg B)$ be $\true$
if and only if $\rm\tau(B)$ is $\false$, and letting $\rm\tau(B\orR C)$ be $\true$ if and only if $\rm\tau(B)$
is $\true$ or $\rm\tau(C)$ is $\true$. A \ii tautology.\/ is a formula A such that $\rm\tau(A)$ is $\true$
for all truth valuations $\tau$ on A; A is a \ii tautological consequence.\/ of $\rm A_1,\ldots,A_\nu$ in
case $\rm A_1\impP\cdots\impP A_\nu\impP A$ is a tautology. Call A and~B \ii tautologically equivalent.\/
in case $\rm A\iffF B$ is a tautology. (Since $\impP$ is associated from right to left,
$\rm A_1\impP\cdots\impP A_\nu\impP A$ is tautologically equivalent to
$\rm A_1\andD\cdots\andD A_\nu\impP A$.)

\everymath={\rm} \everydisplay={\rm}

A \ii numeral.\/ is a term containing no symbols other than S and 0.

Let $\vec x$, y, and z be distinct, let a contain no variables other than those in~$\vec x$, and let b
contain no variables other than those in $\vec x$, y, and z. Then f is \ii defined by primitive
recursion from.\/ a \ii and.\/ b by

\goodbreak

\noin A//1. $f(\vec x,0)=a$

\noin A//2. $f(\vec x,Sy)=b_z\big(f(\vec x,y)\big)$

\smallskip

A \ii construction.\/ of f is a finite sequence $g_1,g_2,\ldots,g_\nu$ where $g_1$ is 0, $g_2$ is S,
f is in the sequence,
and each $g_\mu$ for $3\le\mu\le\nu$ is constructed by primitive recursion from terms containing no
function symbols other than those strictly preceding it in the sequence. The \ii {\rm PR} function symbols.\/ are
those that have a construction.
A \ii {\rm PR} term.\/ is a term in which every function symbol is PR, and a
\ii PR {\rm} number.\/ is a variable-free PR term.

\everymath={} \everydisplay={}

PRA is formulated as a formal system as follows. Its nonlogical symbols are the
PR function symbols. Its nonlogical axioms are the \ii construction axioms.\/ A/1 and A/2
for definitions of PR function symbols, and the \ii successor axioms.


\" \a//3. \neg \'S' x = 0 \"

\" \a//4. \'S' x = \'S' y \imp x = y \"

\smallskip

The logical axioms are \ii reflexivity.\/

\" \a//5. x = x \"

\noin \ii symmetry.\/

\" \a//6. x = y \imp y = x \"

\noin the \ii equality axioms.\/

\noin A//7. $\;\rm x = y \and A_x(x) \imp A_x(y)$

\noin and the \ii propositional axioms.\/

\noin A//8. $\phantom0 \rm A\orR A \impP A$

\noin A//9. $\phantom0 \rm A\impP A\orR B$

\noin A//10. $\rm A\orR B\impP B\orR A$

\noin A//11. $\rm [B\impP C]\impP[A\orR B\impP A\orR C]$

\smallskip

There are three rules of inference:

\smallskip

\quad \ii instance.\/: from A infer an instance of A

\quad \ii modus ponens.\/: from A and $\rm A\impP B$ infer B

\quad \ii induction.\/: from $\rm A_x(0)$ and $\rm A_x(x')\impP A_x(Sx')$ infer $\rm A_x(x)$

\smallskip

As in any formal system,
a \ii proof.\/ is a finite sequence of formulas each of which is either an axiom
or follows from strictly preceding formulas by a rule of inference, and it is a \ii proof of.\/ A
in case A is in the sequence. Decorated $\pi$ stands for a finite sequence of formulas. The notation
$\rm\pi\vdash A$ means that $\pi$ is a proof in PRA of A, while $\rm\vdash A$
means that there is a proof $\pi$ in PRA of A, in which case A is a \ii theorem.\/ of PRA.

\everymath={\rm} \everydisplay={\rm}

The Propositions in this paper are metamathematical in nature; they are statements about PRA whose proofs are
finitary in the strict sense of being expressible in PRA.

\goodbreak

\Pr//  1. {\it Certain familiar devices can be used to extend the notion of proof.

{\rm(i)}~Tautologies are theorems of\/ {\rm PRA}.

{\rm(ii)}~A tautological consequence of theorems of\/ {\rm PRA} is a theorem of\/ {\rm PRA}.

{\rm(iii)}~Previously proved theorems may be cited in proofs.

{\rm(iv)}~Deductions may be used in proofs, as follows. Introduce, in the course of a proof, an arbitrary
formula\/~\ro H, the hypothesis, and follow it by\/ $B_1,\ldots,B_\nu$ where each\/~$B_\mu$ for\/
$1\le\mu\le\nu$ is a theorem
or is a tautological consequence of strictly preceding formulas. Then discharge
the hypothesis by writing\/ \hbox{$H\impP B_\nu$} and never using\/ $H,B_1,\ldots,B_\nu$ again.

{\rm(v)}~Claims may be established, as follows. State a claim\/ \ro A, introduce the hypothesis\/
$\Neg A$, and follow it by\/ $B_1,\ldots,B_\nu$ where these are as in\/~{\rm(iii)} and furthermore a
contradiction is obtained, meaning that for some preceding\/ \ro B is\/ $\Neg B_\nu$.
Then establish the claim by writing\/~\ro A and never using\/ $\Neg A,B_1,\ldots,B_\nu$ again.}

\pf For (i), refer to [HA]. Hilbert and Ackermann give a finitary proof
that in the formal system whose only axioms are
/8--/11 and whose only rule of inference is modus ponens, the theorems are precisely the tautologies.

For (ii), suppose that $\vdash A_1$, \dots, $\vdash A_\nu$ and that $A_1\impP\cdots\impP A_\nu\impP A$
is a tautology, and hence a theorem by (i). Then we have $\vdash A_2\impP\cdots\impP A_\nu\impP A$ by modus
ponens. Proceeding in this way we obtain a proof of A in $\nu$ steps. In other words, tautological
consequence is a derived rule of inference in PRA.

For (iii), just insert the proofs of the cited theorems.

Given a deduction as in (iv), replace $H,B_1,\ldots,B_\nu$ by $H\impP H,H\impP B_1,\ldots,
\hbox{$H\impP B_\nu$}$.
Then $H\impP H$ is a theorem by (i). If $B_\mu$ is a theorem
then $H\impP B_\mu$ is a tautological
consequence of it. If $B_\mu$ is a tautological consequence of strictly preceding formulas, then
$H\impP B_\mu$ is a tautological consequence of them with the $B_\lambda$, for $1\le\lambda<\mu$,
among them replaced by $H\impP B_\lambda$. In this way, by (ii),
we have a proof of $H\impP B_\nu$, proving (iv).
(Notice that no instance of H is taken. This is sometimes expressed by saying that the variables in H
are held constant.)

Given a claim as in (v), proceed as in (iv).
Any formula, in particular~A, is a tautological consequence of $B_\nu$ and $\Neg B_\nu$, so adjoin A
to the deduction. Discharging the hypothesis $\Neg A$ we obtain $\Neg A\impP A$, of which A is a tautological
consequence, proving (v).
(A special case of this is an indirect proof, in which the theorem itself is the claim.) \bul

\everymath={} \everydisplay={}

\smallskip

The only predicate symbol in PRA is =. Nevertheless, we can introduce other predicate symbols as abbreviations.
Given A, let $\rm\vec x$ be its distinct variables in some order, set $\rm p(\vec x)\iffF A$, and let
$\rm p(\vec a)$ abbreviate the instance $\rm A_{\vec x}(\vec a)$. Use decorated p q to stand for predicate
symbols other than =; they occur only in abbreviations.
An \ii explicit definition.\/ of f is $\rm f(\vec x)=c$ where no variable other than those in $\rm\vec x$
occurs in c. Then let $\rm f(\vec a)$ abbreviate $\rm c_{\vec x}(\vec a)$. Explicitly defined function
symbols occur only in abbreviations.

Some formulas are marked $\star$ or ${\star}{\star}$ for emphasis.

\bigskip

\centerline{\bf /   4. Equality}

\medskip

{\it \#2. Construct\/ {\rm Eq} so that\/ $\vdash{\rm Eq}(x,y)=0\iffF x =y$.}

\smallskip

Hilbert and Bernays
construct such a function symbol in \S7 of the first edition of [HB] (1939).
Following a suggestion of Kreisel to establish the basic properties of~$<$ without using addition,
Bernays omitted this explicit construction in the second edition (1968).
We give the surprisingly long construction
of~Eq from the first edition. (Their unary~$\delta$ is our~P, their binary~$\delta$ is our~$-$, and their
$\delta(a,b)+\delta(b,a)$ is our $\ro{Eq}(a,b)$.)

Primitive recursions are labeled with r,
theorems with t, explicit definitions with e, and definitions of PR predicate symbols with~d; if one of these
letters is capitalized, it indicates a schema.
If $\xi$ labels a theorem A, and $\rm x_1$, \dots, $\rm x_\nu$ are the first $\nu$ distinct variables of~A
in the order of
first occurrence, then $\rm\xi;a_1;\ldots;a_\nu$ is the theorem $\rm A_{x_1,\ldots,x_\nu}(a_1,\ldots,a_\nu)$.

\" \ar//12. x + 0 = x \and  x + \'S' y = \'S' ( x + y ) \"

\" \t//13b. 0 + 0 = 0 + 0 \"

\sam13b.
\"\p/13b.
/\'H' \
/5 ; 0 + 0 \
\"

\" \t/13i. x + 0 = 0 + x \imp \'S' x + 0 = 0 + \'S' x \"

\sam13i.
\"\p/13i.
/\'H' : x \
/12 ; \'S' x \
/12 ; 0 ; x \
/12 ; x \
\"

\smallskip

This is an indirect proof. The theorem being proved is a disjunction, $\neg x+0=0+x\orR
\'S' x + 0 = 0 + \'S' x$, and H$:x$ is its negation $x+0=0+x\andD \neg \'S' x + 0 = 0 + \'S' x$
(with $\Neg\Neg$ removed).
The colon indicates that $x$ is to be held fixed with this
introduction of a hypothesis. The remaining formulas in the text proof
together with some equality substitutions (instances of equality axioms /A5) and implicit uses of symmetry a/6,
give a contradiction (marked $QEA$ for quod est absurdum) by tautological consequence, com-\line{pleting
the indirect proof of the theorem. If you are reading this online (it is posted at}\nobreak{\rlap{\hskip5pt\eightit not yet}}
\yy www.math.princeton.edu/$\scriptstyle\sim$nelson/papers/Balrog.pdf\zz),
click on the blue {\nineit Proof\/} link. Otherwise, open a browser to
\yy\pdfklink{www.math.princeton.edu/$\scriptstyle\sim$nelson/proof/}{http://www.math.princeton.edu/~nelson/proof/}\zz\
and click on /13i.pdf.

\" \at/13. x + 0 = 0 + x \"

\smallskip

When t$\xi$ immediately follows t$\xi b$ and t$\xi i$, it is an inference by induction.

\" \t//14b. x + \'S' 0 = \'S' x + 0 \"

\sam14b.
\"\p/14b.
/\'H' : x \
/12 ; x ; 0 \
/12 ; x \
/12 ; \'S' x \
\"

\" \t/14i. x + \'S' y = \'S' x + y \imp x + \'S' \'S' y = \'S' x + \'S' y \"

\sam14i.
\"\p/14i.
/\'H' : x : y \
/12 ; x ; \'S' y \
/12 ; \'S' x ; y \
\"

\" \at/14. x + \'S' y = \'S' x + y \"

\" \t//15b. x + 0 = 0 \imp x = 0 \and 0 = 0 \"

\sam15b.
\"\p/15b.
/\'H' : x \
/12 ; x \
\"

\" \t/15i. [ x + y = 0 \imp x = 0 \and y = 0 ] \imp [ x + \'S' y = 0 \imp x = 0 \and \'S' y = 0 ] \"

\sam15i.
\"\p/15i.
/\'H' : x : y \
/12 ; x ; y \
/3 ; x + y \
\"

\" \at/15. x + y = 0 \imp x = 0 \and y = 0 \"


\" \ar//16. \'P' 0 = 0 \and \'P' \'S' x = x \"


\" \ar//17. x \minus 0 = x \and x \minus \'S' y = \'P' ( x \minus y ) \"

\" \t//18b. \'S' x \minus \'S' 0 = x \minus 0 \"

\sam18b.
\"\p/18b.
/\'H' : x \
/17 ; \'S' x ; 0 \
/17 ; \'S' x \
/16 ; x \
/17 ; x \
\"

\" \t/18i. \'S' x \minus \'S' y = x \minus y \imp \'S' x \minus \'S' \'S' y = x \minus \'S' y \"

\sam18i.
\"\p/18i.
/\'H' : x : y \
/17 ; \'S' x ; \'S' y \
/17 ; x ; y \
\"

\" \at/18. \'S' x \minus \'S' y = x \minus y \"

\" \t//19b. 0 \minus 0 = 0 \"

\sam19b.
\"\p/19b.
/\'H' \
/17 ; 0 \
\"

\" \t/19i. x \minus x = 0 \imp \'S' x \minus \'S' x = 0 \"

\sam19i.
\"\p/19i.
/\'H' : x \
/18 ; x ; x \
\"

\" \at/19. x \minus x = 0 \"

\" \t//20b. \'S' 0 \minus 0 = \'S' 0 \"

\sam20b.
\"\p/20b.
/\'H' \
/17 ; \'S' 0 \
\"

\" \t/20i. \'S' x \minus x = \'S' 0 \imp \'S' \'S' x \minus \'S' x = \'S' 0 \"

\sam20i.
\"\p/20i.
/\'H' : x \
/18 ; \'S' x ; x \
\"

\" \at/20. \'S' x \minus x = \'S' 0 \"

\" \t//21. \'S' x \minus x \ne 0 \"

\sam21.
\"\p/21.
/\'H' : x \
/20 ; x \
/3 ; 0 \
\"

\" \t//22b. 0 \ne 0 \imp 0 = \'S' \'P' 0 \"

\sam22b.
\"\p/22b.
/\'H' \
/5 ; 0 \
\"

\" \t/22i. [ x \ne 0 \imp x = \'S' \'P' x ] \imp [ \'S' x \ne 0 \imp \'S' x = \'S' \'P' \'S' x ] \"

\sam22i.
\"\p/22i.
/\'H' : x \
/16 ; x \
\"

\" \at/22. x \ne 0 \imp x = \'S' \'P' x \"

\" \t//23b. y \minus 0 \ne 0 \imp 0 + ( y \minus 0 ) = y \"

\sam23b.
\"\p/23b.
/\'H' : y \
/17 ; y \
/13 ; y \
/12 ; y \
\"

\" \t/23i. [ y \minus x \ne 0 \imp x + ( y \minus x ) = y ] \imp
[ y \minus \'S' x \ne 0 \imp \'S' x + ( y \minus \'S' x ) = y ] \"

\sam23i.
\"\p/23i.
/\'H' : y : x \
/17 ; y ; x \
/16 ; y \minus x \
/14 ; x ; \'P' ( y \minus x ) \
/22 ; y \minus x \
\"

\" \at/23. y \minus x \ne 0 \imp x + ( y \minus x ) = y \"

\" \t//24. y \minus x = \'S' 0 \imp \'S' x = y \"

\sam24.
\"\p/24.
/\'H' : y : x \
/3 ; 0 \
/23 ; y ; x \
/12 ; x ; 0 \
\"

\" \t//25. x \ne 0 \and \'P' x = 0 \imp x = \'S' 0 \"

\sam25.
\"\p/25.
/\'H' : x \
/22 ; x \
\"

\" \t//26. y \minus x \ne 0 \and y \minus \'S' x = 0 \imp y \minus x = \'S' 0 \"

\sam26.
\"\p/26.
/\'H' : y : x \
/17 ; y ; x \
/25 ; y \minus x \
\"

\" \t//27. y \minus x \ne 0 \imp \'S' x = y \or y \minus \'S' x \ne 0 \"

\sam27.
\"\p/27.
/\'H' : y : x \
/26 ; y ; x \
/24 ; y ; x \
\"

\" \t//28b. x \minus 0 \ne 0 \imp \'S' x \minus 0 \ne 0 \"

\sam28b.
\"\p/28b.
/\'H' : x \
/17 ; \'S' x \
/3 ; x \
\"

\" \t/28i. [ x \minus y \ne 0 \imp \'S' x \minus y \ne 0 ] \imp
[ x \minus \'S' y \ne 0 \imp \'S' x \minus \'S' y \ne 0 ] \"

\sam28i.
\"\p/28i.
/\'H' : x : y \
/18 ; x ; y \
/17 ; x ; y \
/16 \
\"

\" \at/28. x \minus y \ne 0 \imp \'S' x \minus y \ne 0 \"

\" \t//29b. 0 = y \or y \minus 0 \ne 0 \or 0 \minus y \ne 0 \"

\sam29b.
\"\p/29b.
/\'H' : y \
/17 ; y \
\"

\" \t/29i. [ x = y \or y \minus x \ne 0 \or x \minus y \ne 0 ] \imp
[ \'S' x = y \or y \minus \'S' x \ne 0 \or \'S' x \minus y \ne 0 ] \"

\sam29i.
\"\p/29i.
/\'H' : x : y \
/21 ; x \
/27 ; y ; x \
/28 ; x ; y \
\"

\" \at/29. x = y \or y \minus x \ne 0 \or x \minus y \ne 0 \"

\" \t//30. x \minus y = 0 \and y \minus x = 0 \imp x = y \"

\sam30.
\"\p/30.
/\'H' : x : y \
/29 ; x ; y \
\"

\" \de//31. \'Eq' ( x , y ) = ( x \minus y ) + ( y \minus x ) \"

\smallskip

Here at last is the Hilbert-Bernays construction.

\" \t//32. \'Eq' ( x , y ) = 0 \iff x = y \hfill{\star}{\star} \"

\sam32.
\"\p/32.
/\'H' : x : y \
/31 ; x ; y \
/15 ; x \minus y ; y \minus x \
/30 ; x ; y \
/19 ; x \
/12 ; 0 \
/? \'Eq' ( x , y ) \ne 0 \
\"

\smallskip

The ? indicates the introduction of a claim.

\smallskip

{\it \#3. Using the case function symbol\/ \ro C such that\/ $\ro C(x,y,z)$ is y if x is \ro 0 and is z
otherwise, form the characteristic term\/ $\rm\chi A$ so that\/ $\rm\vdash\chi A=0\iffF A$. Construct the
formal system\/ $\chi{\rm PRA}$, equivalent to\/ {\rm PRA}, whose only symbols are variables and\/
{\rm PR} function symbols: the logical connectives and\/ $=$ are eliminated. But continue working in\/ {\rm PRA}.}


\" \ar//33. \'C' ( 0 , y , z ) = y \and \'C' ( \'S' x , y , z ) = z \hfill\star \"

\" \t//34. x \ne 0 \imp \'C' ( x , y , z ) = z \"

\sam34.
\"\p/34.
/\'H' : x : y : z \
/22 ; x \
/33 ; y ; z ; \'P' x \
\"

\" \t//35. \'C' ( x , y , z ) \ne z \imp x = 0 \and \'C' ( x , y , z ) = y \"

\sam35.
\"\p/35.
/\'H' : x : y : z \
/34 ; x ; y ; z \
/33 ; y ; z \
\"

\" \t//36. \'C' ( x , 0 , \'S' 0 ) = 0 \iff x = 0 \"

\sam36.
\"\p/36.
/\'H' : x \
/22 ; x \
/33 ; 0 ; \'S' 0 ; \'P' x \
/3 ; 0 \
/? x \ne 0 \
\"

\" \t//37. \'C' ( x , \'S' 0 , 0 ) = 0 \iff x \ne 0 \"

\sam37.
\"\p/37.
/\'H' : x \
/22 ; x \
/33 ; \'S' 0 ; 0 ; \'P' x \
/3 ; 0 \
/? x = 0 \
\"

\" \t//38. \'C' ( x , 0 , \'S' 0 ) = 0 \or \'C' ( x , 0 , \'S' 0 ) = \'S' 0 \"

\sam38.
\"\p/38.
/\'H' : x \
/22 ; x \
/33 ; 0 ; \'S' 0 ; \'P' x \
\"

\" \t//39. \'C' ( x , \'S' 0 , 0 ) = 0 \or \'C' ( x , \'S' 0 , 0 ) = \'S' 0 \"

\sam39.
\"\p/39.
/\'H' : x \
/22 ; x \
/33 ; \'S' 0 ; 0 ; \'P' x \
\"

\" \t//40. \'C' ( x , 0 , 0 ) = 0 \"

\sam40.
\"\p/40.
/\'H' : x \
/22 ; x \
/33 ; 0 ; 0 ; \'P' x \
\"

\" \t//41a. \'C' \big ( x , 0 , \'C' ( y , 0 , \'S' 0 ) \big ) = 0 \imp x = 0 \or y = 0 \"

\sam41a.
\"\p/41a.
/\'H' : x : y \
/36 ; y \
/22 ; x \
/33 ; 0 ; \'C' ( y , 0 , \'S' 0 ) ; \'P' x \
\"

\" \t/41b. x = 0 \or y = 0 \imp \'C' \big ( x , 0 , \'C' ( y , 0 , \'S' 0 ) \big ) = 0 \"

\sam41b.
\"\p/41b.
/\'H' : x : y \
/36 ; y \
/22 ; x \
/33 ; 0 ; \'C' ( y , 0 , \'S' 0 ) ; \'P' x \
\"

\" \at/41. \'C' \big ( x , 0 , \'C' ( y , 0 , \'S' 0 ) \big ) = 0 \iff x = 0 \or y = 0 \"

\smallskip

{\nineit Proof.} {\ninerm Tautological consequence of /41a and /41b.}

\smallskip

Think of 0 as true and S0 as false.

\" \de//42. x \deq y = \'C' \big ( \'Eq' ( x , y ) , 0 , \'S' 0 \big ) \"

\" \de//43. \dneg x = \'C' ( x , \'S' 0 , 0 ) \"

\" \de//44. x \dor y = \'C' \big ( x , 0 , \'C' ( y , 0 , \'S' 0 ) \big ) \"

\" \t//45. x \deq y = 0 \or x \deq y = \'S' 0 \"

\sam45.
\"\p/45.
/\'H' : x : y \
/42 ; x ; y \
/38 ; \'Eq' ( x , y ) \
\"

\" \t//46. \dneg x = 0 \or \dneg x = \'S' 0 \"

\sam46.
\"\p/46.
/\'H' : x \
/43 ; x \
/39 ; x \
\"

\" \t//47. x \dor y = 0 \or x \dor y = \'S' 0 \"

\sam47.
\"\p/47.
/\'H' : x : y \
/44 ; x ; y \
/38 ; y \
/40 ; x \
/38 ; x \
\"

\" \de//48. x \dand y = \dneg ( \dneg x \dor \dneg y ) \"

\" \de//49. x \dimp y = \dneg x \dor y \"

\" \de//50. x \diff y = ( \dneg x \dor y ) \dand ( x \dor \dneg y ) \"

\smallskip

To each A associate a term $\rm\chi A$, the \ii characteristic term of.\/ A, recursively as follows.
$$\leqalignno{
\rm\chi[a=b]&\;\hbox{ is }\;\rm a\deq b&(// 1)\cr
\rm\chi[\Neg B]&\;\hbox{ is }\;\rm\dneg\chi B&(// 2)\cr
\rm\chi[B\orR C]&\;\hbox{ is }\;\rm\chi B \dor \chi C&(// 3)\cr
}$$
That is, $\rm\chi A$ is obtained by replacing each = by $\deq$, each $\Neg$ by $\dneg$, and each $\orR$
by $\dor$. From this it follows that
$$\rm[\chi A]_{\vec x}(\vec a)\;\ is\;\ \chi[A_{\vec x}(\vec a)]\leqno(// 4)$$

If $\ell$ is a decorated roman letter occurring in an expression schema v, then $\rm v,\ell\slash u$ is the
expression or expression schema obtained by replacing each occurrence of $\ell$ in v by u.

\goodbreak

\Pr//  2. {\it The following are theorem schemata of\/} PRA.

\noin T//51. $\rm\chi A = 0 \iff A$

\noin T//52. $\rm\chi A=0\or\chi A=S0$

\pf We have $\rm /51,A\slash a=b$ by $(/ 1)$ and \rf/42;a;b and $\rm\rf/36;Eq(a,b)$ and $\rm\rf/32;a;b$.
If $\rm /51,A\slash B$ then $\rm/51,A\slash \Neg B$ by~$(/ 2)$ and $\rm\rf/43;\chi B$ and
$\rm\rf/37;\chi B$. If $\rm /51,A\slash B$ and $\rm /51,A\slash C$ then $\rm /51,A\slash B\orR C$
by $(/ 3)$ and $\rm\rf/47;\chi B;\chi C$ and
$\rm\rf/41;\chi B;\chi C$. By metamathematical induction
on the formation of formulas, each $\rm\chi A=0\iffF A$ is a theorem of PRA.

We have $\rm /52,A\slash a=b$ by $(/ 1)$ and \rf/45;a;b; we have $\rm/52,A\slash\Neg B$ by
$(/ 2)$ and $\rm\rf/46;\chi B$; we have \hbox{$\rm/52,A\slash B\orR C$} by $(/ 3)$ and $\rm\rf/47;\chi B;\chi C$.
By definition, every formula A is an equation, negation, or disjunction, so /52 holds. \bul

\goodbreak

\noin T//53. $\rm\chi A=S0\iff\Neg A$

\smallskip

{\nineit Proof.} {\ninerm Tautological consequence of /51 and /52 and \Rf{/3};0.}

\noin T//54. $\rm[A\iffF B]\iff \chi A=\chi B$

\smallskip

{\ninepoint \everymath={}\everymath={\rm}

{\it Proof.} Tautological consequence of /51 and $/51,A\slash B$ and /52 and $/52,A\slash B$ and \Rf{/3};0.

}

\smallskip

Now reformulate PRA as a formal system $\chi$PRA with a simpler data structure.
The symbols of $\chi$PRA are the variables and the PR function symbols. Terms are as before.
A \ii $\chi$-equation.\/ is a term of the form $\rm a\deq b$.
The formulas of $\chi$PRA, called \ii $\chi$-formulas., are
defined recursively as follows. A $\chi$-equation is a $\chi$-formula; if b and c are $\chi$-formulas,
so are $\rm\dneg b$ and $\rm b\dor c$. Decorated $\alpha\ \beta\ \gamma\ \delta$ stand for $\chi$-formulas.
Note that $\chi$ is bijective from formulas of PRA onto formulas of $\chi$PRA; each $\alpha$ is $\rm\chi A$
for a unique A, $\rm\chi^{-1}\alpha$.

Think of the $\chi$-formula $\alpha$ as asserting that the term $\alpha$ is equal to 0.
The axioms of $\chi$PRA are the characteristic terms of the axioms of PRA; the rules of inference of
$\chi$PRA are formed from the rules of inference of PRA by replacing each premise and conclusion by its
characteristic term.

Explicitly, the axioms and rules of inference of $\chi$PRA are as follows, where in A$\chi$/1 and A$\chi$/2,
f is the function symbol defined by A/1 and A/2.

\noin A$\chi$/1.\phantom0  $\rm f(\vec x,y)\deq a$

\noi A$\chi$/2.\phantom0  $\rm f(\vec x,Sy)\deq b_z\big(f(\vec x,y)\big)$

\noi a$\chi$/3.\phantom0\hskip5pt $\dneg\ \ro Sx\deq 0$

\noi a$\chi$/4.\phantom0\hskip5pt $\ro Sx\deq\ro Sy\dimp x\deq y$

\noi a$\chi$/5.\phantom0\hskip5pt $x\deq x$

\noi a$\chi$/6.\phantom0\hskip5pt $x\deq y\dimp y\deq x$

\noi A$\chi$/7.\phantom0  $\rm x\deq y \dand \alpha_x(x)\dimp\alpha_x(y)$

\noi A$\chi$/8.\phantom0  $\alpha\dor\alpha\dimp\alpha$

\noi A$\chi$/9.\phantom0  $\alpha\dimp\alpha\dor\beta$

\noi A$\chi$/10. $\alpha\dor\beta\dimp\beta\dor\alpha$

\noi A$\chi$/11. $(\beta\dimp\gamma)\dimp(\alpha\dor\beta\dimp\alpha\dor\gamma)$

\smallskip

\ii $\chi$-instance.\/: from $\alpha$ infer an instance of $\alpha$

\ii $\chi$-modus ponens.\/: from $\alpha$ and $\alpha\dimp\beta$ infer $\beta$

\ii $\chi$-induction.\/: from $\rm\alpha_x(0)$ and $\rm\alpha_x(x')\dimp\alpha_x(Sx')$ infer $\rm\alpha_x(x)$

\smallskip

A proof in $\chi$PRA is a \ii $\chi$-proof.. Decorated $\sigma$ stands for
a finite sequence of $\chi$-formulas, $\sigma\vdash^\chi\alpha$ asserts that $\sigma$ is a $\chi$-proof
of $\alpha$, and $\vdash^\chi\alpha$ asserts that $\alpha$ is a theorem of $\chi$PRA.

\Pr//  3. {\it
The following are equivalent: $\rm\vdash A$ and\/ $\rm\vdash\chi A=0$ and\/ $\rm\vdash^\chi\chi A$.}

\pf The first two are equivalent by T/51. Let $\rm\pi\vdash A$. Then $\rm\chi\circ\pi\vdash^\chi\chi A$---for
if B in $\pi$ is an axiom, so is $\rm\chi B$, and if B is inferred by a rule of inference, then $\rm\chi B$
is inferred by the corresponding rule. Conversely, if $\rm\sigma\vdash^\chi\chi A$, let $\pi$ consist of
all B of the form $\chi^{-1}\beta=0$ for $\beta$ in $\sigma$. If $\beta$ is an axiom, then $\chi^{-1}\beta$ is
an axiom D of PRA, so B is $\rm\chi D=0$, which is a theorem by $\rm T/51,A\slash D$.
If $\beta$ is inferred by a rule of inference, then B is inferred by the corresponding rule.
Hence $\pi$ is a proof with citation of theorems from the theorem schema T/51, so $\rm\vdash A$. \bul

\smallskip

{\it \#4. Establish primitive recursion by cases, though it will not be used until much later.}

\smallskip

The following Proposition expresses the familiar ``if, else if, \dots, else if, else'' pattern for cases.

\Pr//  4. {\it Let\/ $\rm d$ be\/
$\rm C\chi A_1b_1C\chi A_2b_2\ldots C\chi A_\nu b_\nu c_\nu$ and for\/ $1\le\mu\le\nu$ let\/ $\rm B_\mu$ be\/
$\rm\Neg A_1\andD\cdots\andD\Neg A_\mu$. Then }
$$\eqalign{\rm\vdash[A_1\impP d=b_1]\and&\rm[B_1\andD A_2\impP d=b_2]\and\cdots\and\cr
&\rm[B_{\nu-1}\andD A_\nu\impP d=b_\nu]
\and[B_\nu\impP d=c_\nu]\cr}\leqno(// 5)$$

\pf First let $\nu=1$. Then $\rm d$ is $\rm C\chi A_1b_1c_1$. We have $\rm\chi A_1=0\iffF A_1$ by
$\rm T/51,A\slash A_1$, so $\rm A_1\impP d=b_1$ by $\rm\rf/33;b_1;c_1$. We have $\rm\chi A_1=S0\iffF\Neg A_1$
by $\rm T/53,A\slash A_1$, so $\rm\Neg A_1\impP d=c_1$ by $\rm\rf/33;b_1;c_1;0$, proving the result for $\nu=1$.
Now assume as metamathematical induction hypothesis that the result holds for $\nu-1$ and let $\rm c_1'$ be
$\rm C\chi A_2b_2\ldots C\chi A_\nu b_\nu c_\nu$. By the case $\nu=1$, $\rm A_1\impP d=b_1$ and
$\rm\Neg A_1\impP d=c_1'$, so $(/ 5)$ holds by the induction hypothesis. \bul

\smallskip

If f is defined by the primitive recursion $\rm f(\vec x,y)=a\andD f(\vec x,Sy)=d$ (where d is
as in the Proposition) then we say that $(/ 5)$ holds by \ii primitive recursion by cases..

\bigskip

\centerline{\bf /   5. Bounded quantifiers}

\medskip

{\it \#5. Given a formula\/ \ro A and a variable\/ \ro x, construct the\/ {\rm PR} function symbol\/
$\rm \mu\down A$ so that\/ $\rm\mu\down A(x)$ (with the other variables in\/ \ro A not indicated in the
notation) finds the first\/ $\rm x'$, if any, such that\/ $\rm A_x(x')$ holds.}

\smallskip

Introduce the PR predicate symbols $\le$ (\ii less than.\/) and $\lt$ (\ii strictly less than.\/).

\" \d//55. x \le y \iff x \minus y = 0 \"

\" \d//56. x \lt y \iff x \le y \and x \ne y \"

\" \t//57. x \le 0 \imp x = 0 \"

\sam57.
\" \p/57.
/\'H' : x \
/55\fw ; x ; 0 \
/17 ; x \
\"

\smallskip

If d$\xi$ is $\rm p(\vec x)\iffF D$, then $\xi\fw$ is $\rm\Neg p(\vec x)\orR D$ (which is
$\rm p(\vec x)\impP D$) and $\xi\bw$ is $\rm p(\vec x)\orR \Neg D$ (which is tautologically equivalent
to $\rm D\impP p(\vec x)$).

\" \t//58. 0 \le x \"

\sam58.
\"\p/58.
/\'H' : x \
/55\bw ; 0 ; x \
/23 ; 0 ; x \
/15 ; x ; 0 \minus x \
\"

\" \t//59. \neg \'S' x \le x \"

\sam59.
\"\p/59.
/\'H' : x \
/55\fw ; \'S' x ; x \
/21 ; x \
\"

\" \t//60. x \le x \"

\sam60.
\"\p/60.
/\'H' : x \
/55\bw ; x ; x \
/19 ; x \
\"

\" \t//61. x \le y \imp x \lt y \or x = y \"

\sam61.
\"\p/61.
/\'H' : x : y \
/56\bw ; x ; y \
\"

\" \t//62. x \lt y \or x = y \imp x \le y \"

\sam62.
\"\p/62.
/\'H' : x : y \
/56\fw ; x ; y \
/60 ; x \
\"

\" \t//63. x \le \'S' x \"

\sam63.
\"\p/63.
/\'H' : x \
/17 ; x ; x \
/19 ; x \
/16 \
/55\bw ; x ; \'S' x \
\"

\" \t//64b. \'P' 0 \le 0 \"

\sam64b.
\"\p/64b.
/\'H' \
/16 \
/60 ; 0 \
\"

\" \t/64i. \'P' x \le x \imp \'P' \'S' x \le \'S' x \"

\sam64i.
\"\p/64i.
/\'H' : x \
/16 ; x \
/63 ; x \
\"

\" \at/64. \'P' x \le x \"

\" \t//65. x \lt y \imp y \minus x \ne 0 \"

\sam65.
\"\p/65.
/\'H' : x : y \
/56\fw ; x ; y \
/55\fw ; x ; y \
/29 ; x ; y \
\"

\" \t//66. x \lt y \imp y = x + ( y \minus x ) \"

\sam66.
\"\p/66.
/\'H' : x : y \
/65 ; x ; y  \
/23 ; y ; x \
\"

\" \t//67. x = y \imp y = x + ( y \minus x ) \"

\sam67.
\"\p/67.
/\'H' : x : y \
/19 ; y \
/12 ; x \
\"

\" \t//68. x \le y \imp y = x + ( y \minus x ) \"

\sam68.
\"\p/68.
/\'H' : x : y \
/56\bw ; x ; y \
/67 ; x ; y \
/66 ; x ; y \
\"

\" \t//69b. x \minus ( y + 0 ) = ( x \minus y ) \minus 0 \"

\sam69b.
\"\p/69b.
/\'H' : x : y \
/12 ; y \
/17 ; x \minus y \
\"

\" \t/69i. [ x \minus ( y + z ) = ( x \minus y ) \minus z ]
\imp [ x \minus ( y + \'S' z ) = ( x \minus y ) \minus \'S' z ] \"

\sam69i.
\"\p/69i.
/\'H' : x : y : z \
/12 ; y ; z \
/17 ; x ; y + z \
/17 ; x \minus y ; z \
\"

\" \at/69. x \minus ( y + z ) = ( x \minus y ) \minus z \"

\" \t//70b. 0 \minus 0 = 0 \"

\sam70b.
\"\p/70b.
/\'H' \
/19 ; 0 \
\"

\" \t/70i. 0 \minus x = 0 \imp 0 \minus \'S' x = 0 \"

\sam70i.
\"\p/70i.
/\'H' : x \
/17 ; 0 ; x \
/16 \
\"

\" \at/70. 0 \minus x = 0 \"

\" \t//71. x \le x + y \"

\sam71.
\"\p/71.
/\'H' : x : y \
/69 ; x ; x ; y \
/19 ; x \
/70 ; y \
/55\bw ; x ; x + y \
\"

\" \t//72b. x + ( y + 0 ) = ( x + y ) + 0 \"

\sam72b.
\"\p/72b.
/\'H' : x : y \
/12 ; y \
/12 ; x + y
\"

\" \t/72i. x + ( y + z ) = ( x + y ) + z
\imp x + ( y + \'S' z ) = ( x + y ) + \'S' z \"

\sam72i.
\"\p/72i.
/\'H' : x : y :z \
/12 ; y ; z \
/12 ; x ; y + z \
/12 ; x + y ; z \
/4 ; x + ( y + z ) ; ( x + y ) + z \
\"

\" \at/72. x + ( y + z ) = ( x + y ) + z \"

\" \t//73. x \le y \and y \le z \imp x \le z \"

\sam73.
\"\p/73.
/\'H' ; x : y : z \
/68 ; x ; y \
/68 ; y ; z \
/72 ; x ; y \minus x ; z \minus y \
/71 ; x ; ( y \minus x ) + ( z \minus y ) \
\"

\" \t//74. x \lt y \imp y \minus x \ne 0 \"

\sam74.
\"\p/74.
/\'H' : x : y \
/56\fw ; x ; y \
/68 ; x ; y \
/12 ; x \
\"

\" \t//75. y \minus x \ne 0 \imp x \lt y \"

\sam75.
\"\p/75.
/\'H' : y : x \
/23 ; y ; x \
/71 ; x ; y \minus x \
/56\bw ; x ; y \
/19 ; x \
\"

\" \t//76. x \le y \and y \le x \imp x = y \"

\sam76.
\"\p/76.
/\'H' : x : y \
/55\fw ; x ; y \
/55\fw ; y ; x \
/30 ; x ; y \
\"

\" \t//77. x \lt y \or x = y \or y \lt x \"

\sam77.
\"\p/77.
/\'H' : x : y \
/75 ; x ; y \
/75 ; y ; x \
/30 ; x ; y \
\"

\" \t//78. x \lt y \imp \neg y \le x \"

\sam78.
\"\p/78.
/\'H' : x : y \
/56\fw ; x ; y \
/76 ; x ; y \
\"

\" \t//79. \neg y \le x \imp x \lt y \"

\sam79.
\"\p/79.
/\'H' : y : x \
/77 ; x ; y \
/56\fw ; y ; x \
/60 ; x \
\"

\" \t//80. \neg [ x \lt y \and y \le x ] \"

\sam80.
\"\p/80.
/\'H' : x : y \
/78 ; x ; y \
\"

\" \t//81. \'S' x \le y \imp x \lt y \"

\sam81.
\"\p/81.
/\'H' : x : y \
/56\bw ; x ; y \
/63 ; x \
/73 ; x ; \'S' x ; y \
/55\fw ; \'S' x ; x \
/21 ; x
\"

\smallskip

\everymath={\rm} \everydisplay={\rm}

Let $\vec y$ be the distinct variables of A other than x in the order of first occurrence,
and let $\mu\down A(a)$ abbreviate $\mu\down A(\vec y,a)$.
Define the PR function symbol $\rm\mu\down A$ by

\noin R//82. $\mu\down A(0)=C(\chi[A_x(0)],0,S0) \and \hfil\break
\mu\down A(Sx)=C\big(\chi[\mu\down A(x)\le x],\mu\down A(x),C(\chi[A_x(Sx)],Sx,SSx)\big)\hfill\star $

\smallskip

Remark that if $x_1$ is any variable other than those in $\vec y$ and $A_1$ is $A_x(x_1)$, then
$\mu_{A_1}$ is the same as $\mu\down A$.
We have $\mu\down A(x)=Sx$ until an $x'$ (if any) is found such that $A_x(x')$ holds,
after which it remains~$x'$ forever, as we now demonstrate (with implicit uses of equality axioms and symmetry).

\noin $T//83b$. $\mu\down A(0)\le0\imp A_x\big(\mu\down A(0)\big)$

\smallskip

{\ninepoint \everymath={}\everymath={\rm}

{\it Proof.} Suppose not. Then

.1\quad $\mu\down A(0)\le0$

.2\quad $\Neg A_x\big(\mu\down A(0)\big)$

\noi By .1 and $\rf/57;\mu\down A(0)$,

.3\quad $\mu\down A(0)=0$

\noi By .3 and /82,

.4\quad $C\big(\chi[A_x(0)],0,S0\big)=0$

\noi By .4 and $\rf/36;\chi[A_x(0)]$,

.5\quad $\chi[A_x(0)]=0$

\noi By .5 and $/51,A\slash A_x(0)$,

.6\quad $A_x(0)$

\noi By .2 and .3,

.7\quad $\Neg A_x(0)$

\noi QEA by .6 and .7.

}

\goodbreak

\noin $T/83i$.  $\big[\mu\down A(x)\le x \impP A_x\big(\mu\down A(x)\big)\big]\imp
\big[\mu\down A(Sx)\le Sx\impP A_x\big(\mu\down A(Sx)\big)\big]$

\smallskip

{\ninepoint \everymath={}\everymath={\rm}

{\it Proof.} Suppose not. Then

.1\quad $\mu\down A(x)\le x \impP A_x\big(\mu\down A(x)\big)$

.2\quad $\mu\down A(Sx)\le Sx$

.3\quad $\Neg A_x\big(\mu\down A(Sx)\big)$

Claim: $\neg\mu\down A(x)\le x$. Suppose not. Then

\qquad .4\quad $\mu\down A(x)\le x$

\noi By .4 and $/51,A\slash\mu\down A(x)\le x$,

\qquad .5\quad $\chi[\mu\down A(x)]=0$

\noi By .4 and /82,

\qquad .6\quad $\mu\down A(Sx)=C\big(0,\mu\down A(x),C(\chi[A_x(Sx)],Sx,SSx)\big)$

\noi By .6 and $\rf/33;\mu\down A(x)$,

\qquad .7\quad $\mu\down A(Sx)=\mu\down A(x)$

\noi By .1 and .4,

\qquad .8\quad $A_x\big(\mu\down A(x)\big)$

\noi By .7 and .8,

\qquad .9\quad $A_x\big(\mu\down A(Sx)\big)$

\noi The claim is proved by .3 and .9, and .4--.9 will not be used again.

.10\quad $\neg\mu\down A(x)\le x$

\noi By .10 and $/53,A\slash \mu\down A(x)\le x$,

.11\quad $\chi[\mu\down A(x)\le x]=S0$

\noi By .11 and /82,

.12\quad $\mu\down A(Sx)=C\big(S0,\mu\down A(c),C(\chi[A_x(Sx)],Sx,SSx)\big)$

\noi By .12 and $\rf/33;\mu\down A(x);C\big(\chi[A_x(Sx)],Sx,SSx\big);0$,

.13\quad $\mu\down A(Sx)=C\big(\chi[A_x(Sx)],Sx,SSx\big)$

\noi By .2 and \rf/59;Sx,

.14\quad $\mu\down A(Sx)\ne SSx$

\noi By .14 and $\rf/35;\chi[A_x(Sx)];Sx;SSx$,

.15\quad $\chi[A_x(Sx)]=0\andD C\big(\chi[A_x(Sx)],Sx,SSx\big)=Sx$

\noi By .15 and $/51,A\slash A_x(Sx)$,

.16\quad $A_x(Sx)$

\noi By .13 and .15,

.17\quad $\mu\down A(Sx)=Sx$

\noi By .16 and .17,

.18\quad $A_x\big(\mu\down A(Sx)\big)$

\noi QEA by .3 and .18.

}

\noin $T/83.$ $\mu\down A(x)\le x \imp A_x\big(\mu\down A(x)\big)$

\noin $T//84{\mit b}.$ $A_x(0) \imp \mu\down A(0)\le0$

\smallskip

{\ninepoint \everymath={}\everymath={\rm}

{\it Proof.} Suppose not. Then

.1\quad $A_x(0)$

.2\quad $\neg \mu\down A(0)\le0$

\noi By .1 and $/49,A\slash A_x(0)$,

.3\quad $\chi[A_x(0)=0]=0$

\noi By .3 and /82,

.4\quad $\mu\down A(0)=C(0,S0,SS0)$

\noi By .4 and \Rf{/40};0,

.5\quad $\mu\down A(0)=0$

\noi By .5 and \Rf{/60};0,

.6\quad $\mu\down A(0)\le0$

\noi QEA by .2 and .7.

}

\noin $T/84{\mit i}$. $[A_x(x)\impP \mu\down A(x)\le x]\imp[A_x(Sx)\impP\mu\down A(Sx)\le Sx]$

\smallskip

{\ninepoint \everymath={}\everymath={\rm}

{\it Proof.} Suppose not. Then (the induction hypothesis is not needed in this proof)

.1\quad $A_x(Sx)$

.2\quad $\neg \mu\down A(Sx)\le Sx$

\noi By .1 and $/51,A\slash A_x(Sx)$,

.3\quad $\chi[A_x(Sx)]=0$

\noi  Claim: $\mu\down A(x)\le x$. Suppose not. Then

\qquad .4\quad $\neg \mu\down A(x)\le x$

\noi By .4 and $/53,A\slash\mu\down x(x)\le x$,

\qquad .5\quad $\chi[\mu\down A(x)\le x]=S0$

\noi By .5 and /82,

\qquad .6\quad $\mu\down A(Sx)=C\big(S0,\mu\down A(x),C(\chi[A_x(Sx)],Sx,SSx)\big)$

\noi By .6 and $\rf/33;\mu\down A(x);C(\chi[A_x(Sx)],Sx,SSx);0$,

\qquad .7\quad $\mu\down A(Sx)=C\big(\chi[A_x(Sx)],Sx,SSx\big)$

\noi By .3 and .7,

\qquad .8\quad $\mu\down A(Sx)=C(0,Sx,SSx)$

\noi By .8 and \rf/33;Sx;SSx,

\qquad .9\quad $\mu\down A(Sx)=Sx$

\noi By .9 and \rf/60;Sx,

\qquad .10\quad $\mu\down A(Sx)\le Sx$

\noi The claim is proved by .2 and .10, and .4--.10 will not be used again.

.11\quad $\mu\down A(x)\le x$

\noi By .11 and $/51,A\slash \mu\down A(x)\le x$,

.12\quad $\chi[\mu\down A(x)\le x]=0$

\noi By .12 and /82,

.13\quad $\mu\down A(Sx)= C\big(0,\mu\down A(x),C(\chi[A_x(Sx)],Sx,SSx)\big)$

\noi By .13 and $\rf/33;\mu\down A(x);C\big(\chi[A_x(Sx)],Sx,SSx\big)$,

.14\quad $\mu\down A(Sx)=\mu\down A(x)$

\noi By .14 and .11,

.15\quad $\mu\down A(Sx)\le x$

\noi By .15 and \rf/63;x and $\rf/73;\mu\down A(x);x;Sx$,

.16\quad $\mu\down A(Sx)\le Sx$

\noi QEA by .2 and .16.

}

\goodbreak
\noin T/84. $A_x(x) \imp \mu\down A(x) \le x$

\noin T//85. $A_x(x)\imp A_x\big(\mu\down A(x)\big) \and \mu\down A(x) \le x$

\smallskip

{\nineit Proof.} {\ninerm  By /84 and /83.}

\goodbreak

\noin $T//86{\mit b}$. $\mu\down A(x)\le x \imp \mu\down A(x+0)=\mu\down A(x)$

\smallskip

{\ninepoint \everymath={}\everymath={\rm}

{\it Proof.} By \rf/12;x and $\rf/5;\mu\down A(x)$.

}

\noin $T/86{\mit i}$. $[\mu\down A(x)\le x \impP \mu\down A(x+w)=\mu\down A(x)] \imp
[\mu\down A(x)\le x \impP \mu\down A(x+Sw)=\mu\down A(x)]$

\smallskip

{\ninepoint \everymath={}\everymath={\rm}

{\it Proof.} Suppose not. Then

.1\quad $\mu\down A(x)\le x\impP \mu\down A(x+w)=\mu\down A(x)$

.2\quad $\mu\down A(x)\le x$

.3\quad $\mu\down A(x+Sw)\ne \mu\down A(x)$

\noi By .2 and .1,

.4\quad $\mu\down A(x+w)=\mu\down A(x)$

\noi By \rf/12;x;w,

.5\quad $x+Sw=S(x+w)$

\noi By .4 and .2,

.6\quad $\mu\down A(x+w)\le x$

\noi By .6 and \rf/71;x;w and $\rf/73;\mu\down A(x+w);x;x+w$,

.7\quad $\mu\down A(x+w)\le x+w$

\noi By .7 and $/51,A\slash \mu\down A(x+w)\le x+w$,

.8\quad $\chi[\mu\down A(x+w)\le x+w]=0$

\noi By .8 and $\rf/82;x+w$,

.9\quad $\mu\down A\big(S(x+w)\big)=C\big(0,\mu\down A(x+w),C(\chi[A_x(S(x+w))],S(x+w),SS(x+w))\big)$

\noi By .9 and .5 and $\rf/33;\mu\down A(x+w);C(\chi[A_x(S(x+w))],S(x+w),SS(x+w))$,

.10\quad $\mu\down A(x+Sw)=\mu\down A(x+w)$

\noi By .10 and .4,

.11\quad $\mu\down A(x+Sw)=\mu\down A(x)$

\noi QEA by .3 and .11.

}

\noin T/86. $\mu\down A(x)\le x \imp \mu\down A(x+w)=\mu\down A(x)$

\noin T//87. $\mu\down A(x)\le x \and x\le z \imp \mu\down A(z)=\mu\down A(x)$

\smallskip

{\ninepoint \everymath={}\everymath={\rm}

{\it Proof.} Suppose not. Then

.1\quad $\mu\down A(x)\le x$

.2\quad $x\le z$

.3\quad $\mu\down A(z)\ne\mu\down A(x)$

\noi By .2 and \rf/68;x;z,

.4\quad $z=x+(z-x)$

\noi By .1 and .4 and $/86;x;z-x$,

.5\quad $\mu\down A(z)=\mu\down A(x)$

\noi QEA by .3 and .5.

}

\smallskip

{\it \#6. Introduce bounded quantifiers\/ $\rm\exists x\leE b\,A$ and\/ $\rm\forall x\leE b\,A$, where\/ \ro x
does not occur in\/ \ro b, as instances of\/ \ro A, using\/ $\rm\mu\down A$.}

\smallskip

If x does occur in b, let $\exists x\leE b\,A$ abbreviate $A_x(\mu\down A(b))$,
and let $\forall x\leE b\,A$ abbreviate $A_x(\mu\down{\Neg A}(b))$. Then
$\forall x\leE b\,A$ is tautologically equivalent to $\Neg\exists x\leE b\,\Neg A$, since the latter is the
double negation of the former.
Call $\exists x\leE b$ and $\forall x\le b$ \ii bounded quantifiers.\/; they occur only in abbreviations,
and x does not occur in $\exists x\leE b\,A$ or $\forall x\leE b\,A$. Each is an instance of A.

\noin T//88. $x\le b \and A_x(x) \imp \exists x\leE b\,A $

\smallskip

{\ninepoint \everymath={}\everymath={\rm}

{\it Proof.} Suppose not. Then

.1\quad $x\le b$

.2\quad $A_x(x)$

.3\quad $\neg A_x\big(\mu\down A(b)\big)$

\noi By .2 and /85,

.4\quad $A_x\big(\mu\down A(x)\big) \andD \mu\down A(x) \le x $

\noi By .4 and .1 and /87;x;b,

.5\quad $\mu\down A(b)=\mu\down A(x)$

\noi By .5 and .4,

.6\quad $A_x\big(\mu\down A(b)\big)$

\noi QEA by .3 and .6.

}

\smallbreak

For the converse we have a specific number, $\mu\down A(b)$, that is less than b and satisfies~A.

\noin T//89. $\exists x\leE b\,A\imp\mu\down A(b)\le b \and A_x\big(\mu\down A(b)\big) \hfill\star $

\smallskip

{\ninepoint \everymath={}\everymath={\rm}

{\it Proof.} By $/85;\mu\down A(b)$ (recalling that $\exists x\leE b\,A$ is $A_x(\mu\down A(b))$).

}

\noin T//90. $\forall x\leE b\,A \imp[x\le b\imp A_x(x)] \hfill\star $

\smallskip

{\ninepoint \everymath={}\everymath={\rm}

{\it Proof.} Tautologically equivalent to $/88,A\slash\Neg A$.

}

\noin T//91. $\exists x\leE c\,A \and c \le b \imp \exists x\leE b\,A$

\smallskip

{\ninepoint \everymath={}\everymath={\rm}

{\it Proof.} Suppose not. Then

.1\quad $A_x\big(\mu\down A(c)\big)$

.2\quad $c\le b$

.3\quad $\neg A_x\big(\mu\down A(b)\big)$

\noi By $/89,b\slash c$,

.4\quad $\mu\down A(c)\le c \andD A_x\big(\mu\down A(b)\big)$

\noi By $/87;\mu\down A(c);c;b$ and .2,

.5\quad $\mu\down A(b)=\mu\down A(c)$

\noi QEA by .1 and .3 and .5.

}

\noin T//92. $\forall x\leE b\,A \and c\le b \imp \forall x\leE c\,A$

\smallskip

{\ninepoint \everymath={}\everymath={\rm}

{\it Proof.} Tautologically equivalent to $/91,A\slash\Neg A$.

}

\smallskip

{\it \#7. Replace induction as a rule of inference by an axiom schema, and use this to construct
formal systems\/ {\rm PRA*} and $\chi${\rm PRA*}, equivalent to\/ {\rm PRA}, in which the only rule of
inference is, respectively, modus ponens or $\chi$-modus ponens, and such that variable-free theorems have
variable-free proofs. But continue working in\/ {\rm PRA}; $\chi${\rm PRA*} will be arithmetized later.}

\noin $T//93\mit b$. $A_x(0) \and \forall x'\leE Px[A_x(x')\impP A_x(Sx')]\imp A_x(0)$

{\nineit Proof.} {\ninerm Tautology.}

\noin $T/93\mit i$. $\{A_x(0) \and \forall x'\leE Px[A_x(x')\impP A_x(Sx')]\imp A_x(x)\}\imp \hfill\break
\{A_x(0) \and \forall x'\leE PSx[A_x(x')\impP A_x(Sx')]\imp A_x(Sx)\}$

\smallskip

{\ninepoint \everymath={}\everymath={\rm}

{\it Proof.} Suppose not. Then

.1\quad $A_x(0)\andD\forall x'\leE Px[A_x(x')\impP A_x(Sx')]\impP A_x(x)$

.2\quad $A_x(0)$

.3\quad $\forall x'\leE PSx[A_x(x)\impP A_x(Sx')]$

.4\quad $\neg A_x(Sx)$

\noi By .2 and .1,

.5\quad $\forall x'\leE Px[A_x(x')\impP A_x(Sx')]\impP A_x(x)$

\noi By \rf/16;x,

.6\quad $PSx=x$

\noi By .6 and .3,

.7\quad $\forall x'\leE x[A_x(x')\impP A_x(Sx')]$

\noi By \rf/60;x,

.8\quad $x\le x$

\noi By .8 and .7 and $/90,x\slash x',b\slash x,A\slash A_x(x)\impP A_x(Sx)$,

.9\quad $A_x(x)\impP A_x(Sx)$

\noi By \rf/64;x,

.10\quad $Px\le x$

\noi By .10 and .7 and $/90,x\slash x',b\slash Px,A\slash A_x(x)\impP A_x(Sx)$,

.11\quad $\forall x'\leE Px[A_x(')\impP A_x(Sx')]$

\noi By .11 and .5,

.12 $A_x(x)$

\noi By .12 and .9,

.13\quad $A_x(Sx)$

\noi QEA by .4 and .13.

}

\noin T/93. $A_x(0) \and \forall x'\leE Px[A_x(x') \impP A_x(Sx')] \imp A_x(x) \hfill{\star}{\star}$

\smallskip

Now simplify PRA even further, via formal theories PRA* and $\chi$PRA*. The symbols, terms, and formulas of
PRA* are those of PRA; its axioms are all instances of axioms of PRA and of T/93; its only rule of inference
is modus ponens. The symbols, terms, and formulas of $\chi$PRA* are those of $\chi$PRA; its axioms are the
characteristic terms of the axioms of PRA*; its only rule of inference is $\chi$-modus ponens.
Proofs and theorems of these two systems are indicated by $\vdash^*$ and $\vdash^{\chi*}$.

\Pr//  5. {\it {\rm(i)} $\vdash A$ if and only if\/ $\vdash^* A$.

{\rm(ii)} A variable-free theorem of\/ {\rm PRA*} has a variable-free proof.

{\rm(iii)} $\vdash^*A$ if and only if\/ $\vdash^{\chi*}\chi A$.

{\rm(vi)} A variable-free theorem of\/ $\rm\chi PRA^*$ has a variable-free proof.}

\pf For (i), let $\pi\vdash A$ and let $B_x(x)$ be the first formula in $\pi$ that is inferred by induction,
from $B_x(0)$ and $B_x(x')\impP B_x(Sx')$. Then we have proofs of these two premises without using induction.
By an instance of the second premise we have \hbox{$\forall x\leE Px[B_x(x')\impP B_x(Sx')]$}, so we have
$B_x(0)\andD \forall x\leE Px[B_x(x')\impP B_x(Sx')]$ by tautological consequence. By $/93,A\slash B$ we have
$B_x(x)$. Proceeding in this way, by metamathematical induction on the number of inferences by induction
in~$\pi$, we have
a proof~$\pi_1$ of~A from the axioms of PRA* without using induction. Now consider the first formula
$C_{\vec x}(\vec a)$ of $\pi_1$ that is inferred by instance, from C. Let $\pi_2$ the the part of $\pi_1$
strictly preceding $C_{\vec x}(\vec a)$.
Let $D'$ be $D_{\vec x}(\vec a)$, and let $\pi_2'$ consist of all $D'$ for D in $\pi_2$.
Then $\pi_2\pi_2'C'$ is a proof from the axioms of PRA* of $C'$
(i.e., $C_{\vec x}(\vec a)$) without using induction
or instance, because if D is an axiom of PRA* so is $D'$ (since an instance of an instance is an instance),
and if D is inferred from $D_1$ and \hbox{$D_1\impP D$} by modus ponens, then $D'$ is inferred from $D_1'$ and
$D_1'\impP D'$ by modus ponens (since \hbox{$[D_1\impP D]'$} is $D_1'\impP D'$).
(Don't delete $\pi_2$, because later
formulas in $\pi_1$ may be inferred by instance from a formula in it.) Proceeding in this way. by
metamathematical induction on the number of inferences by instance in $\pi_1$, we obtain a proof of A in PRA*.
The converse direction of (i) is trivial, since the axioms of PRA* are theorems of PRA and the rule of inference
of PRA* is a rule of inference of PRA.

For (ii), let A be variable-free with $\pi\vdash^*A$. Let $D^\circ$ be the formula obtained by replacing all
variables in D by 0, and let $\pi^\circ$ consist of all $D^\circ$ for D in $\pi$. Then $\pi^\circ$ is a
variable-free proof in PRA* of A, because if D is an axiom of PRA* so is $D^\circ$; if D is inferred
by modus ponens from $D_1$ and $D_1\impP D$, then $D^\circ$ is inferred by modus ponens from $D_1^\circ$
and $D_1^\circ\impP D^\circ$; and $A^\circ$ is A.

For (iii), let $\pi\vdash^* A$. Then $\chi\circ\pi\vdash^{\chi*}\chi A$---for
if B in $\pi$ is an axiom, so is $\chi B$, and if B is inferred by modus ponens, then $\chi B$
is inferred by $\chi$-modus ponens. Conversely, if $\sigma\vdash^{\chi*}\chi A$, let $\pi$ consist of
all B of the form $\chi^{-1}\beta=0$ for $\beta$ in $\sigma$. If $\beta$ is an axiom, then $\chi^{-1}\beta$ is
an axiom D of PRA*, which is a theorem of PRA,
so B is $\chi D=0$, which is a theorem of PRA by $T/51,A\slash D$.
If $\beta$ is inferred by $\chi$-modus ponens, then B is inferred by modus ponens.
Hence $\pi$ is a proof in PRA with citation of theorems, so $\vdash A$.
Consequently, $\vdash^*A$ by (i).

For (iv), let $\chi A$ be a variable-free theorem of $\chi$-PRA. By (iii), A is a variable-free theorem of
PRA*, which by (ii) has a variable-free proof $\pi$ in PRA*. As shown in the proof of (iii),
$\chi\circ\pi\vdash^{\chi*}\chi A$, and this proof is variable-free. \bul

\smallskip

{\it \#8. Construct the least number principle as a derived rule of inference.}

\smallskip

If the recursion for $\mu\down A$ finds an $x'$ such that $A_x(x')$ holds, then $x'$ is the least number for
which $A_x(x')$ holds. This leads to the least number principle, a powerful form of induction.

\goodbreak

\noin T//94. $A_x(x) \imp [y\lt\mu\down A(x) \imp \Neg A_x(y)]$

\smallskip

{\ninepoint \everymath={}\everymath={\rm}

{\it Proof.} Suppose not. Then

\quad .1 $A_x(x)$

\quad .2 $y\lt\mu\down A(x)$

\quad .3 $A_x(y)$

\noi By .3 and /85;y,

\quad .4 $A_x\big(\mu\down A(y)\big)$

\quad .5 $\mu\down A(y)\le y$

\noi By .1 and /85;x,

\quad .6 $A_x\big(\mu\down A(x)\big)$

\quad .7 $\mu\down A(x)\le x$

\noi By .2 and $/56\fw;y;\mu\down A(x)$,

\quad .8 $y\le\mu\down A(x)$

\noi By .8 and .7 and $\rf/73;\mu\down A(x);x$,

\quad .9 $y\le x$

\noi By .5 and .9 and /87;y;x,

\quad .10 $\mu\down A(x)=\mu\down A(y)$

\noi By .2 and .5 and .10,

\quad .11 $y\lt\mu\down A(y)\andD \mu\down A(y)\le y$

\noi QEA by .11 and $\rf/80;y;\mu\down A(y)$.}

\smallskip

Let $\forall x\ltT b\,A$ (where x does not occur in b) abbreviate $\forall x\leE b[x\lt b\impP A]$.
Then $\forall x\ltT b\,A$ is an instance of $x\lt b\impP A$. The following theorem schema follows from T/94
by instance.

\noin T//95. $A_x(x)\imp\forall y\ltT \mu\down A(x)[\Neg A_x(y)]$

\noin T//96. $A_x(x) \imp \mu\down A(\mu\down A(x))=\mu\down A(x)$

\smallskip

{\ninepoint \everymath={}\everymath={\rm}

{\it Proof.} Suppose not. Then

\quad .1 $A_x(x)$

\quad .2 $\mu\down A(\mu\down A(x))\ne\mu\down A(x)$

\noi By .1 and /85,

\quad .3 $A_x(\mu\down A(x))$

\noi By .3 and $/84;\mu\down A(x)$,

\quad .4 $A_x(\mu\down A(\mu\down A(x)))$

\quad .5 $\mu\down A(\mu\down A(x))\le\mu\down A(x)$

\noi By .2 and .5 and $/56\bw;\mu\down A(x);\mu\down A(\mu\down A(x))$,

\quad .6 $\mu\down A(\mu\down A(x))\lt\mu\down A(x)$

\noi By .3 and .6 and $\rf/94;\mu\down A(x);\mu\down A(\mu\down A(x))$,

\quad .7 $\Neg A_x(\mu\down A(\mu\down A(x)))$

\noi QEA by .4 and .7.}

\Pr//  6. {\it If\/ \ro A does not have a least counterexample, then\/ \ro A holds. That is,
if\/ $\vdash\Neg\{\Neg A_x(z)\andD\forall y\ltT z[A_x(y)]\}$ then\/ $\vdash A_x(x)$.}

\pf Suppose .0 $\vdash\Neg\{\Neg A_x(z)\andD\forall y\ltT z[A_x(y)]\}$.
Then we prove $A_x(x)$ as follows. Suppose not. Then

\goodbreak

\quad .1 $\Neg A_x(x)$

\noi By .1 and $/85,A\slash\Neg A$,

\quad .2 $\Neg A(\mu\down A(x))$

By .2 and $/95,A\slash\Neg A,x\slash\mu\down A(x)$,

\quad .3 $\forall y\ltT\mu\down{\Neg A}(\mu\down{\Neg A}(x))[A_x(y)]$

\noi By .3 and .1 and $/96,A\slash\Neg A$,

\quad .4 $\forall y\ltT\mu\down{\Neg A}(x)[A_x(y)]$

\noi By $.0;\mu\down{\Neg A}(x)$,

\quad .5 $\Neg\{\Neg A_x(\mu\down{\Neg A}(x))\andD \forall y\ltT\mu\down{\Neg A}(x)[A_x(y)]\}$

\noi By .5 and .2,

\quad .6 $\Neg\forall y\ltT\mu\down{\Neg A}(x)[A_x(y)]$

\noi QEA by .4 and .6. \bul

\smallskip

The derived rule of inference, from $\Neg\{\Neg A_x(z)\andD
\forall y\ltT z[A_x(y)]\}$ infer A, is the \ii least number principle..

\smallskip

{\it \#9. Construct definition of function symbols with uniqueness condition and bounded existence condition.}

\Pr//  7. {\it Let\/ \ro A contain no variables other than the distinct variables\/ $\vec y$ and\/~\ro x.
The uniqueness condition\/ {\rm(UC)} is\/ $A_x(x)\andD A_x(x')\impP x =x'$, where\/ $x'$ is distinct from the\/
$\vec y$ and\/ \ro x. The existence condition\/ {\rm(EC)} is\/ $\exists x\leE b\,A$, where\/ \ro b contains no
variables other than those in\/ $\vec y$. If\/ {\rm UC} and\/ {\rm EC} are theorems, let\/ $f(\vec y)$,
where \ro f is a new function symbol, abbreviate\/ $\mu\down A(b)$. Then $\vdash f(\vec y)=x\iffF A$.}

\pf By definition, EC is $\mu\down A(b)$, and since it is a theorem, $\vdash A_x\big(\mu\down A(b)\big)$---i.e.,
\hbox{$\vdash A_x\big(f(\vec y)\big)$}---by T/89. Consequently, $\vdash f(\vec y)=x\impP A$. The converse holds
by UC. \bul

\smallskip

Such function symbols are \ii defined function symbols.; they occur only in abbreviations.

\bigbreak

\centerline{\bf /   6. Strings}

\medskip

\everymath={} \everydisplay={}

A string is a concatenation of bits. Identify the number $x$ with
the string consisting of the ones and zeros following the initial one in the binary
representation of $\ro Sx$. We implement this in PRA.

\smallskip

{\it \#10. Prove that every non-zero number\/ $\ro Sx$ can be written uniquely as\/ $\ro Qx+\ro Rx$
where\/ \ro Qx is a power of two and\/ \ro Rx is strictly less than \ro Qx.}

\" \t//97b. 0 + 0 = 0 \"

\sam97b.
\"\p/97b.
/\'H' \
/12 ; 0 \
\"

\" \t/97i. 0 + x = x \imp 0 + \'S' x = \'S' x \"

\sam97i.
\"\p/97i.
/\'H' : x \
/12 ; 0 ; x \
\"

\" \at/97. 0 + x = x \"

\" \t//98b. x + 0 = 0 + x \"

\sam98b.
\"\p/98b.
/\'H' : x \
/12 ; x \
/97 ; x \
\"

\" \t/98i. x + y = y + x \imp x + \'S' y = \'S' y + x \"

\sam98i.
\"\p/98i.
/\'H' : x : y \
/12 ; x ; y \
/12 ; y ; x \
/14 ; y ; x \
\"

\" \at/98. x + y = y + x \"


\smallskip

Introduce \ii multiplication.. As usual, $\cdot$ binds more tightly than $+$.

\" \ar//99. x \cdot 0 = 0 \and x \cdot \'S' y = x + x \cdot y \"

\" \t//100. x \cdot 0 = 0 \and x \cdot \'S' y = x \cdot y + x \"

\sam100.
\"\p/100.
/\'H' : x : y \
/99 ; x ; y \
/98 ; x ; x \cdot y \
\"

\" \t//101b. x \cdot ( y + 0 ) = x \cdot y + x \cdot 0 \"

\sam101b.
\"\p/101b.
/\'H' : x : y \
/12 ; y \
/100 ; x \
/12 ; x \cdot y \
\"

\" \t/101i. x \cdot ( y + z ) = x \cdot y + x \cdot z \imp
x \cdot ( y + \'S' z ) = x \cdot y + x \cdot \'S' z \"

\sam101i.
\"\p/101i.
/\'H' : x : y : z \
/12 ; y ; z \
/100 ; x ; z \
/100 ; x ; y + z \
/72 ; x \cdot y ; x \cdot z ; x \
\"

\" \at/101. x \cdot ( y + z ) = x \cdot y + x \cdot z \"

\" \t//102b. x \cdot ( y \cdot 0 ) = ( x \cdot y ) \cdot 0 \"

\sam102b.
\"\p/102b.
/\'H' : x : y \
/100 ; y \
/100 ; x \
/100 ; x \cdot y \
\"

\" \t/102i. x \cdot ( y \cdot z ) = ( x \cdot y ) \cdot z \imp
x \cdot ( y \cdot \'S' z ) = ( x \cdot y ) \cdot \'S' z \"

\sam102i.
\"\p/102i.
/\'H' : x : y : z \
/100 ; y ; z \
/100 ; x \cdot y ; z \
/101 ; x ; y \cdot z ; y \
\"

\" \at/102. x \cdot ( y \cdot z ) = ( x \cdot y ) \cdot z \"

\" \t//103b. 0 \cdot 0 = 0 \"

\sam103b.
\"\p/103b.
/\'H' \
/100 ; 0 \
\"

\" \t/103i. 0 \cdot x = 0 \imp 0 \cdot \'S' x = 0 \"

\sam103i.
\"\p/103i.
/\'H' : x \
/100 ; 0 ; x \
/12 ; 0 \
\"

\" \at/103. 0 \cdot x = 0 \"

\" \t//104b. \'S' x \cdot 0 = x \cdot 0 + 0 \"

\sam104b.
\"\p/104b.
/\'H' : x \
/100 ; \'S' x \
/100 ; x \
/12 ; 0 \
\"

\" \t/104i. \'S' x \cdot y = x \cdot y + y \imp \'S' x \cdot \'S' y = x \cdot \'S' y + \'S' y \"

\sam104i.
\"\p/104i.
/\'H' : x : y \
/100 ; x ; y \
/100 ; \'S' x ; y \
/72 ; x \cdot y ; x ; \'S' y \
/72 ; x \cdot y ; y ; \'S' x \
/14 ; y ; x \
/98 ; \'S' y ; x \
\"

\" \at/104. \'S' x \cdot y = x \cdot y + y \"

\" \t//105b. x \cdot 0 = 0 \cdot x \"

\sam105b.
\"\p/105b.
/\'H' : x \
/100 ; x \
/103 ; x \
\"

\" \t/105i. x \cdot y = y \cdot x \imp x \cdot \'S' y = \'S' y \cdot x \"

\sam105i.
\"\p/105i.
/\'H' : x : y \
/100 ; x ; y \
/104 ; y ; x \
\"

\" \at/105. x \cdot y = y \cdot x \"

\" \t//106. ( x + y ) \cdot z = x \cdot z + y \cdot z \"

\sam106.
\"\p/106.
/\'H' : x : y : z \
/105 ; x + y ; z \
/101 ; z ; x ; y \
/105 ; z ; x \
/105 ; z ; y \
\"

\" \t//107. x \le y \imp x \cdot z \le y \cdot z \"

\sam107.
\"\p/107.
/\'H' : x : y : z \
/68 ; x ; y \
/106 ; x ; y \minus x ; z \
/71 ; x \cdot z ; ( y \minus x ) \cdot z \
\"

\" \t//108. x \lt y \imp x + ( y \minus x ) = y \and y \minus x \ne 0 \"

\sam108.
\"\p/108.
/\'H' : x : y \
/68 ; x ; y \
/56\fw ; x ; y \
/12 ; x \
\"

\" \t//109. x \lt \'S' y \imp x \le y \"

\sam109.
\"\p/109.
/\'H' : x : y \
/108 ; x ; \'S' y \
/22 ; \'S' y \minus x \
/12 ; x ; \'P' ( \'S' y \minus x ) \
/4 ; x + \'P' ( \'S' y \minus x ) ; y \
/71 ; x ; \'P' ( \'S' y \minus x ) \
\"

\"\t//110. x \le \'S' y \imp x \le y \or x = \'S' y \"

\sam110.
\"\p/110.
/\'H' : x : y \
/61 ; x ; \'S' y \
/109 ; x ; y \
\"

\" \t//111. \neg [ x \lt y \and y \lt \'S' x ] \"

\sam111.
\"\p/111.
/\'H' : x : y \
/109 ; y ; x \
/78 ; x ; y \
\"

\" \t//112. x \le y \imp \'S' x \le \'S' y \"

\sam112.
\"\p/112.
/\'H' : x : y \
/68 ; x ; y \
/12 ; x ; y \minus x \
/14 ; x ; y \minus x \
/71 ; \'S' x ; y \minus x \
\"

\" \t//113. x \lt y \imp \'S' x \lt \'S' y \"

\sam113.
\"\p/113.
/\'H' : x : y \
/112 ; x ; y \
/56\fw ; x ; y \
/56\bw ; \'S' x ; \'S' y \
/4 ; x ; y \
\"

\" \t//114. x \lt y \imp \'S' x \le y \"

\sam114.
\"\p/114.
/\'H' : x : y \
/108 ; x ; y \
/22 ; y \minus x \
/14 ; x ; \'P' ( y \minus x ) \
/71 ; \'S' x ; \'P' ( y \minus x ) \
\"

\" \de//115. 1 = \'S' 0 \"

\" \de//116. 2 = \'S' \'S' 0 \"

\" \t//117. 1 \cdot x = x \"

\sam117.
\"\p/117.
/\'H' : x \
/115 \
/104 ; 0 ; x \
/97 ; x \
/103 ; x \
\"

\" \t//118. 2 \cdot x = x + x \"

\sam118.
\"\p/118.
/\'H' : x \
/116 \
/115 \
/104 ; 1 ; x \
/117 ; x \
\"

\" \t//119b. y + 0 = z + 0 \imp y = z \"

\sam119b.
\"\p/119b.
/\'H' : y : z \
/12 ; y \
/12 ; z \
\"

\" \t/119i. [ y + x = z + x \impP y = z ] \imp [ y + \'S' x = z + \'S' x \impP y = z ] \"

\sam119i.
\"\p/119i.
/\'H' : y : x : z \
/12 ; y ; x \
/12 ; z ; x \
/4 ; y + x ; z + x \
\"

\" \at/119. y + x = z + x \imp y = z \"

\" \t//120. x + y = x + z \imp y = z \"

\sam120.
\"\p/120.
/\'H' : x : y : z \
/98 ; x ; y \
/98 ; x ; z \
/119 ; y ; x ; z \
\"

\" \t//121. y \ne 0 \imp x \lt x + y \"

\sam121.
\"\p/121.
/\'H' : y : x \
/71 ; x ; y \
/56\bw ; x ; x + y \
/120 ; x ; 0 ; y \
/12 ; x \
\"

\" \t//122. x \le y \and y \lt z \imp x \lt z \"

\sam122.
\"\p/122.
/\'H' : x : y : z \
/56\fw ; y ; z \
/73 ; x ; y ; z \
/56\bw ; x ; z \
/80 ; y ; x \
\"

\" \t//123. x \lt y \imp y \ne 0 \"

\sam123.
\"\p/123.
/\'H' : x : y \
/56\fw ; x ; y \
/57 ; x \
\"

\" \t//124. x \lt y \imp \'S' x \lt 2 \cdot y \"

\sam124.
\"\p/124.
/\'H' : x : y \
/118 ; y \
/114 ; x ; y \
/123 ; x ; y \
/121 ; y ; y \
/122 ; \'S' x ; y ; y + y \
\"

\" \t//125. x \lt \'S' x \"

\sam125.
\"\p/125.
/\'H' : x \
/59 ; x \
/79 ; \'S' x ; x \
\"


\goodbreak

\smallskip Introduce \ii exponentiation..

\" \ar//126. x \up 0 = 1 \and x \up \'S' y = x \cdot ( x \up y ) \"

\" \t//127b. 0 \lt 2 \up 0 \"

\sam127b.
\"\p/127b.
/\'H' \
/116 \
/115 \
/126 ; 2 \
/125 ; 0 \
\"

\" \t/127i. x \lt 2 \up x \imp \'S' x \lt 2 \up \'S' x \"

\sam127i.
\"\p/127i.
/\'H' : x \
/126 ; 2 ; x \
/124 ; x ; 2 \up x \
\"

\" \at/127. x \lt 2 \up x \"

\" \t//128. x \le 2 \up x \"

\sam128.
\"\p/128.
/\'H' : x \
/127 ; x \
/56\fw ; x ; 2 \up x \
\"

\" \d//129. q \' is a power of two' \iff \exists x \leE q [ 2 \up x = q ] \hfill\star \"

\smallskip

Although ``is a power of two'' contains English words, it is a formal predicate symbol, written in
suffix notation. The negation of ``$q \' is a power of two'$'' is ``$\neg q \' is a power of two'$'',
whereas ``$q \' is not a power of two'$'' is not even an expression.

\" \t//130. 1 \' is a power of two' \"

\sam130.
\"\p/130.
/\'H' \
/115 \
/129\bw ; 1 ; 0 \
/128 ; 0 \
/126 ; 2 \
\"

\smallskip

Here is the explanation of the notation $\rm/129\bw;1;0$ in the proof. The formula $/129\bw$ is

\smallskip

\qquad $q \' is a power of two' \or \Neg\exists x\leE q[2\up x = q]$

\noin which is tautologically equivalent to

\smallskip

\qquad $q \' is a power of two' \or \forall x\leE q[2\up x\ne q]$

\noin Then $\rm/129\bw;1$ is tautologically equivalent to

\smallskip

\qquad $1 \' is a power of two' \or \forall x\leE 1[2\up x\ne1]$

\noin This is a variable-free formula and the ;0 in $\rm/129\bw;1;0$ is an implicit use of T/90.
That is, $\rm/129\bw;1;0$ is

\smallskip

\qquad $1 \' is a power of two' \or [0\le1\imp 2\up0\ne1]$

\smallskip

In general, $\rm\forall x\leE b\,A$ obeys the expected rule: we can substitute any term c, indicated
by ;c, and obtain $\rm c\le b \impP A_x(c)$.

\" \t//131. 2 \up x \' is a power of two' \"

\sam131.
\"\p/131.
/\'H' : x \
/116 \
/129\bw ; 2 \up x ; x \
/128 ; x \
/5 ; 2 \up x \
\"

\" \t//132. x \cdot y = 0 \imp x = 0 \or y = 0 \"

\sam132.
\"\p/132.
/\'H' : x : y \
/22 ; y \
/99 ; x ; \'P' y \
/15 ; x ; x \cdot \'P' y \
\"

\" \t//133b. 2 \up 0 \ne 0 \"

\sam133b.
\"\p/133b.
/\'H' \
/116 \
/115 \
/126 ; 2 \
/3 ; 1 \
\"

\" \t/133i. 2 \up x \ne 0 \imp 2 \up \'S' x \ne 0 \"

\sam133i.
\"\p/133i.
/\'H' : x \
/116 \
/115 \
/126 ; 2 ; x \
/132 ; 2 ; 2 \up x \
/3 ; 1 \
\"

\" \at/133. 2 \up x \ne 0 \"

\" \t//134. \neg 0 \' is a power of two' \"

\sam134.
\"\p/134.
/\'H' \
/129\fw ; 0 : x \
/133 ; x \
\"

\smallskip

Here is the explanation of the notation $/129\fw;0:x$ in the proof. The formula $/129\fw$ is

\smallskip

\qquad $\neg q \' is a power of two' \or \exists x\leE q[2\up x =q]$

\noin so $/129\fw;0$ is

\smallskip

\qquad $\neg 0 \' is a power of two' \or \exists x\leE 0[2\up x =0]$

\noin Let $x$ abbreviate $\mu^{\phantom'}_{[{\rm SS} 0\up x =0]}(0)$ and implicitly use T/89.
Then $/129\fw;0:x$ is

\smallskip

\qquad $\neg 0 \' is a power of two' \or [x\leE 0\and2\up x =0]$

\smallskip

In general, $\rm\exists x\leE b\,A$ obeys the expected rule: we can choose any y not previously used in
the proof, hold it fixed, indicated by :y, and obtain $\rm y\le b \andD A_x(y)$.

\" \t//135. q \' is a power of two' \imp 2 \cdot q \' is a power of two' \"

\sam135.
\"\p/135.
/\'H' : q \
/129\fw ; q : x \
/129\bw ; 2 \cdot q ; \'S' x \
/126 ; 2 ; x \
/128 ; \'S' x \
\"

\" \t//136b. x \up ( y + 0 ) = ( x \up y ) \cdot ( x \up 0 ) \"

\sam136b.
\"\p/136b.
/\'H' : x : y \
/12 ; y \
/126 ; x \
/117 ; x \up y \
/105 ; x \up y ; 1 \
\"

\" \t/136i. x \up ( y + z ) = ( x \up y ) \cdot ( x \up z ) \imp
x \up ( y + \'S' z ) = ( x \up y ) \cdot ( x \up \'S' z ) \"

\sam136i.
\"\p/136i.
/\'H' : x : y : z \
/12 ; y ; z \
/126 ; x ; y + z \
/126 ; x ; z \
/102 ; x ; x \up y ; x \up z \
/102 ; x \up y ; x ; x \up z \
/105 ; x ; x \up y \
\"

\" \at/136.  x \up ( y + z ) = ( x \up y ) \cdot ( x \up z ) \"

\goodbreak

\" \t//137. x \ne 0 \imp x \lt 2 \cdot x \"

\sam137.
\"\p/137.
/\'H' : x \
/22 ; x \
/125 ; \'P' x \
/124 ; \'P' x ; x \
\"

\" \t//138. 0 \le x \"

\sam138.
\"\p/138.
/\'H' : x \
/97 ; x \
/71 ; 0 ; x \
\"

\" \t//139. x \ne 0 \imp 1 \le x \"

\sam139.
\"\p/139.
/\'H' : x \
/116 \
/115 \
/22 ; x \
/138 ; \'P' x \
/112 ; 0 ; \'P' x \
\"

\" \t//140. x \ne 0 \imp 2 \le 2 \up x \"

\sam140.
\"\p/140.
/\'H' : x \
/116 \
/115 \
/127 ; x \
/114 ; x ; 2 \up x \
/139 ; x \
/112 ; 1 ; x \
/73 ; 2 ; \'S' x ; 2 \up x \
\"

\" \t//141a. x \cdot z = y \cdot z \and x \lt y \and z \ne 0 \imp x = y \"

\sam141a.
\"\p/141a.
/\'H' : x : z : y \
/108 ; x ; y \
/106 ; x ; y \minus x ; z \
/132 ; y \minus x ; z \
/12 ; x \cdot z \
/120 ; x \cdot z ; 0 ; ( y \minus x ) \cdot z \
\"

\" \t/141. x \cdot z = y \cdot z \and z \ne 0 \imp x = y \"

\sam141.
\"\p/141.
/\'H' ; x ; z ; y \
/141a ; x ; z ; y \
/141a ; y ; z ; x \
/77 ; x ; y \
\"

\" \t//142. x \le y \imp z \cdot x \le z \cdot y \"

\sam142.
\"\p/142.
/\'H' :x : y : z \
/107 ; x ; y ; z \
/105 ; z ; x \
/105 ; z ; y \
\"

\" \t//143. x \ne 0 \imp x \lt x \cdot 2 \"

\sam143.
\"\p/143.
/\'H' : x \
/137 ; x \
/105 ; x ; 2 \
\"

\" \t//144. x \lt y \and y \le z \imp x \lt z \"

\sam144.
\"\p/144.
/\'H' : x : y : z \
/56\fw ; x ; y \
/73 ; x ; y ; z \
/56\bw ; x ; z \
/80 ; x ; y \
\"

\" \t//145. x \lt y \imp 2 \up x \lt 2 \up y \"

\sam145.
\"\p/145.
/\'H' : x : y \
/108 ; x ; y \
/136 ; 2 ; x ; y \minus x \
/140 ; y \minus x \
/142 ; 2 ; 2 \up ( y \minus x ) ; 2 \up x \
/133 ; x \
/143 ; 2 \up x \
/144 ; 2 \up x ; ( 2 \up x ) \cdot 2 ; 2 \up y \
\"

\"\t//146. x \le y \imp 2 \up x \le 2 \up y \"

\sam146.
\"\p/146.
/\'H' : x : y \
/61 ; x ; y \
/145 ; x ; y \
/56\fw ; 2 \up x ; 2 \up y \
/60 ; 2 \up x \
\"

\" \t//147. 2 \up x \lt 2 \up y \imp x \lt y \"

\sam147.
\"\p/147.
/\'H' : x : y \
/77 ; x ; y \
/56\fw ; 2 \up x ; 2 \up y \
/145 ; y ; x \
/80 ; 2 \up x ; 2 \up y \
/56\fw ; 2 \up y ; 2 \up x \
\"

\" \t//148a. 2 \up x = 2 \up y \imp \neg x \lt y \"

\sam148a.
\"\p/148a.
/\'H' : x : y \
/145 ; x ; y \
/56\fw ; 2 \up x ; 2 \up y \
\"

\" \t/148. 2 \up x = 2 \up y \imp x = y \"

\sam148.
\"\p/148.
/\'H' : x : y \
/148a ; x ; y \
/148a ; y ; x \
/77 ; x ; y \
\"

\" \t//149. \neg [ q \' is a power of two' \and q' \' is a power of two' \and q \lt q' \and q' \lt
2 \cdot q ] \"

\sam149.
\"\p/149.
/\'H' : q : q' \
/135 ; q \
/129\fw ; q : x \
/129\fw ; q' : y \
/129\fw ; 2 \cdot q : z \
/126 ; 2 ; x \
/147 ; x ; y \
/147 ; y ; z \
/148 ; \'S' x ; z \
/111 ; x ; y \
\"

\" \d//150. \'p'_{/150} ( q , x ) \iff q \' is a power of two' \and q \le \'S' x \and
\'S' x \lt 2 \cdot q \"

\" \t//151a. \'p'_{/150} ( q , x ) \and \'p'_{/150} ( q' , x ) \imp \neg q \lt q' \"

\sam151a.
\"\p/151a.
/\'H' : q : x : q' \
/150\fw ; q ; x \
/150\fw ; q' ; x \
/149 ; q ; q' \
/122 ; q' ; \'S' x ; 2 \cdot q \
/56\fw ; q ; q' \
\"

\" \t/151. \'p'_{/150} ( q , x ) \and \'p'_{/150} ( q' , x ) \imp q = q' \"

\sam151.
\"\p/151.
/\'H' : q : x : q' \
/151a ; q ; x ; q' \
/151a ; q' ; x ; q \
/77 ; q ; q' \
\"

\" \t//152. \'p'_{/150} ( q , x ) \and \'S' \'S' x \lt 2 \cdot q \imp \'p'_{/150} ( q , \'S' x ) \"

\sam152.
\"\p/152.
/\'H' : q : x \
/150\fw ; q ; x \
/150\bw ; q ; \'S' x \
/63 ; \'S' x \
/73 ; q ; \'S' x ; \'S' \'S' x \
\"

\" \t//153. \'p'_{/150} ( q , x ) \and \neg \'S' \'S' x \lt 2 \cdot q \imp
\'p'_{/150} ( 2 \cdot q , \'S' x ) \"

\sam153.
\"\p/153.
/\'H' : q : x \
/150\fw ; q ; x \
/150\bw ; 2 \cdot q ; \'S' x \
/135 ; q \
/124 ; \'S' x ; 2 \cdot q \
/114 ; \'S' x ; 2 \cdot q \
/56\bw ; \'S' \'S' x ; 2 \cdot q \
/60 ; \'S' \'S' x \
\"

\" \t//154. x \le y \imp x \le \'S' y \"

\sam154.
\"\p/154.
/\'H' : x : y \
/63 ; y \
/73 ; x ; y ; \'S' y \
\"

\" \t//155. \'p'_{/150} ( q , x ) \imp \exists q_1 \leE \'S' \'S' x [ \'p'_{/150} ( q_1 , \'S' x ) ] \"

\sam155.
\"\p/155.
/\'H' : q : x ^A \
/^A ; q \
/^A ; 2 \cdot q \
/152 ; q ; x \
/153 ; q ; x \
/150\fw ; q ; x \
/150\fw ; 2 \cdot q ; \'S' x \
/154 ; q ; \'S' x \
\"

\smallskip

The formula $\ro H:q:x$ is tautologically equivalent to $\'p'_{/150}(q,x)\andD\forall q_1\leE\'S' \'S' x
[\Neg\'p'_{/150}(q_1,\'S'x)]$. The superscript $^A$ labels the remnant
$\forall q_1\leE\'S' \'S' x[\Neg\'p'_{/150}(q_1,\'S'x)]$ for later substitution of values for $q_1$.

\" \t//156b. \exists q \leE 1 [ \'p'_{/150} ( q , 0 ) ] \"

\sam156b.
\"\p/156b.
/\'H' ; 1 \
/116 \
/115 \
/150\bw ; 1 ; 0 \
/130 \
/60 ; 1 \
/3 ; 0 \
/137 ; 1 \
\"

\" \t/156i. \exists q \leE \'S' x [ \'p'_{/150} ( q , x ) ] \imp
\exists q_1 \leE \'S' \'S' x [ \'p'_{/150} ( q_1 , \'S' x ) ] \"

\sam156i.
\"\p/156i.
/\'H' : x : q ^A \
/155 ; q ; x : q_1 \
/^A ; q_1 \
\"

\smallskip

Note that $\exists q \leE \'S' x [ \'p'_{/150} ( q , x ) ]$ is $\mu\down{\ro A}(\'S' x)$ where A is
$\'p'_{/150} ( q , x )$, and $\exists q_1 \leE \'S' \'S' x [ \'p'_{/150} ( q_1 , \'S' x ) ]$ is
$\mu\down{\ro A_1}(\'S' x)$ where $\rm A_1$ is $\ro A_q(q_1)$. By the remark after R/82 they are the same
formula, so induction applies.

\" \at/156. \exists q \leE \'S' x [ \'p'_{/150} ( q , x ) ] \"

\" \d//157. \'Q' x = q \iff \'p'_{/150} ( q , x ) \"

\smallskip

{\ninerm By /151 (UC) and /156 (EC) (and Proposition/  7).}

\" \t//158. \'Q' x \' is a power of two' \and \'Q' x \le \'S' x \and \'S' x \lt 2 \cdot \'Q' x
\hfill\star \"

\sam158.
\"\p/158.
/\'H' : x \
/157\fw ; x ; \'Q' x \
/150\fw ; \'Q' x ; x \
/5 ; \'Q' x \
\"

\smallskip

$\'Q'x$ expresses $\lfloor\log_2(x+1)\rfloor$, the largest power of two less than $x+1$.

\" \t//159. q \' is a power of two' \and q \le \'S' x \and \'S' x \lt 2 \cdot q \imp q = \'Q' x
\hfill\star \"

\sam159.
\"\p/159.
/\'H' : q : x \
/150\bw ; q ; x \
/157\bw ; x ; q \
\"

\" \de//160. \'R' x = \'S' x \minus \'Q' x \"

\" \t//161. \'S' x = \'Q' x + \'R' x \"

\sam161.
\"\p/161.
/\'H' : x \
/158 ; x \
/160 ; x \
/68 ; \'Q' x ; \'S' x \
\"

\" \t//162. x \lt y \imp \exists w \leE y [ x + w = y \and w \ne 0 ] \"

\sam162.
\"\p/162.
/\'H' : x : y ; y \minus x \
/108 ; x ; y \
/98 ; x ; y \minus x \
/71 ; y \minus x ; x \
\"

\" \t//163. x + y = z \and y \ne 0 \imp x \lt z \"

\sam163.
\"\p/163.
/\'H' : x : y : z \
/56\bw ; x ; z \
/71 ; x ; y \
/12 ; x \
/120 ; x ; y ; 0 \
\"

\" \t//164. x + y \lt x + z \imp y \lt z \"

\sam164.
\"\p/164.
/\'H' : x : y : z \
/162 ; x + y ; x + z : w \
/72 ; x ; y ; w \
/120 ; x ; y + w ; z \
/163 ; y ; w ; z \
\"

\" \t//165. \'R' x \lt \'Q' x \"

\sam165.
\"\p/165.
/\'H' : x \
/161 ; x \
/158 ; x \
/118 ; \'Q' x \
/164 ; \'Q' x ; \'R' x ; \'Q' x \
\"

\" \t//166. \'S' x = \'Q' x + \'R' x \and \'Q' x \' is a power of two' \and \'R' x \lt \'Q' x \hfill\star \"

\sam166.
\"\p/166.
/\'H' : x \
/161 ; x \
/158 ; x \
/165 ; x \
\"

\" \t//167. x \le y \imp \exists w \leE y [ x + w = y ] \"

\sam167.
\"\p/167.
/\'H' : x : y ; y \minus x \
/68 ; x ; y \
/98 ; x ; y \minus x \
/71 ; y \minus x ; x \
\"

\" \t//168. x + y \le x + z \imp y \le z \"

\sam168.
\"\p/168.
/\'H' : x : y : z \
/167 ; x + y ; x + z : w \
/72 ; x ; y ; w \
/120 ; x ; y + w ; z \
/71 ; y ; w \
\"

\" \t//169. x + y \le x + z \imp y \le z \"

\sam169.
\"\p/169.
/\'H' : x : y : z \
/167 ; x + y ; x + z : w \
/72 ; x ; y ; w \
/120 ; x ; y + w ; z \
/71 ; y ; w \
\"

\" \t//170. x \lt y \imp x + z \lt y + z \"

\sam170.
\"\p/170.
/\'H' : x : y : z \
/79 ; y + z ; x + z \
/98 ; y ; z \
/98 ; x ; z \
/168 ; z ; y ; x \
/80 ; x ; y \
\"

\" \t//171. \'S' x = q + r \and q \' is a power of two' \and r \lt q \imp q = \'Q' x \and r = \'R' x
\hfill\star \"

\sam171.
\"\p/171.
/\'H' : x : q : r \
/71 ; q ; r \
/118 ; q \
/98 ; q ; r \
/170 ; r ; q ; q \
/159 ; q ; x \
/120 ; \'Q' x ; r ; \'R' x \
/161 ; x \
\"

\smallskip

{\it \#11. Introduce concatenation of strings and prove a few of its properties.}


\" \de//172. x \oplus y = \'P' ( \'S' x \cdot \'Q' y + \'R' y ) \hfill\star \"

\smallskip

Why is $\oplus$ called concatenation? As an example, consider the string 101 (i.e., the number $x$ such that
S$x$ in binary is 1101) and the string 01 (i.e., the number $y$ such that S$y$ in binary is 101). In binary,
$\'Q'y$ is 100, $\'R'y$ is 1, and $\'S'x\cdot\'Q'y+\'R'y$ is 110101, so $x\oplus y$ is the string
10101 (i.e., the number whose successor in binary is 110101).

String arithmetic is analogous to number arithmetic, with one zero, $\ep$, but with two successors:
$x\mapsto x\oplus \u0$ and $x\mapsto x\oplus\u1$. Concatenation is the string analogue of addition.
We have founded string arithmetic on number arithmetic, but we need to develop it to the point that
it becomes independent of this foundation.

\" \t//173. x \ne 0 \and y \lt z \imp y \lt x \cdot z \"

\sam173.
\"\p/173.
/\'H' : x : y : z \
/22 ; x \
/104 ; \'P' x ; z \
/98 ; \'P' x \cdot z ; z \
/71 ; z ; \'P' x \cdot z \
/144 ; y ; z ; x \cdot z \
\"

\" \t//174. x \lt y \imp z + x \lt z + y \"

\sam174.
\"\p/174.
/\'H' : x : y : z \
/98 ; z ; x \
/98 ; z ; y \
/170 ; x ; y ; z \
\"

\" \t//175. r \lt q \and r' \lt q' \imp r \cdot q' + r' \lt q \cdot q' \"

\sam175.
\"\p/175.
/\'H' : r : q : r' : q' \
/162 ; r ; q : w \
/162 ; r' ; q' : w' \
/101 ; r ; r' ; w' \
/106 ; r ; w ; r' + w' \
/173 ; w ; r' ; q' \
/174 ; r' ; w \cdot q' ; r \cdot r' + r \cdot w' \
\"

\" \t//176. \'R' x \cdot \'Q' y + \'R' y \lt \'Q' x \cdot \'Q' y \"

\sam176.
\"\p/176.
/\'H' : x : y \
/166 ; x \
/166 ; y \
/175 ; \'R' x ; \'Q' x ; \'R' y ; \'Q' y \
\"

\" \t//177. q_1 \' is a power of two' \and q_2 \' is a power of two' \imp q_1 \cdot q_2 \' is a power of two' \"

\sam177.
\"\p/177.
/\'H' : q_1 : q_2 \
/129\fw ; q_1 : x_1 \
/129\fw ; q_2 : x_2 \
/136 ; 2 ; x_1 ; x_2 \
/131 ; x_1 + x_2 \
\"

\" \t//178. \'Q' x \ne 0 \"

\sam178.
\"\p/178.
/\'H' : x \
/158 ; x \
/134 \
\"

\" \t//179. \'S' ( x \oplus y ) = \'S' x \cdot \'Q' y + \'R' y \"

\sam179.
\"\p/179.
/\'H' ; x ; y \
/172 ; x ; y \
/22 ; \'S' x \cdot \'Q' y + \'R' y \
/15 ; \'S' x \cdot \'Q' y ; \'R' y \
/3 ; x \
/178 ; y \
/132 ; \'S' x ; \'Q' y \
\"

\" \t//180. \'Q' ( x \oplus y ) = \'Q' x \cdot \'Q' y \and \'R' ( x \oplus y ) = \'R' x \cdot \'Q' y + \'R' y \"

\sam180.
\"\p/180.
/\'H' : x : y \
/179 ; x ; y \
/166 ; x \
/166 ; y \
/106 ; \'Q' x ; \'R' x ; \'Q' y \
/72 ; \'Q' x \cdot \'Q' y ; \'R' x \cdot \'Q' y ; \'R' y \
/177 ; \'Q' x ; \'Q' y \
/176 ; x ; y \
/171 ; \break x \oplus y ; \'Q' x \cdot \'Q' y ; \'R' x \cdot \'Q' y + \'R' y \
\"

\" \t//181. \'Q' \big ( x \oplus ( y \oplus z ) \big ) = \'Q' \big ( ( x \oplus y ) \oplus z \big ) \"

\sam181.
\" \p/181.
/\'H' : x : y : z \
/180 ; y ; z \
/180 ; x ; y \
/180 ; x ; y \oplus z \
/180 ; x \oplus y ; z \
/102 ; \'Q' x ; \'Q' y ; \'Q' z \
\"

\" \t//182. \'R' \big ( x \oplus ( y \oplus z ) \big ) = \'R' \big ( ( x \oplus y ) \oplus z \big ) \"

\sam182.
\" \p/182.
/\'H' : x : y : z \
/180 ; y ; z \
/180 ; x ; y \
/180 ; x ; y \oplus z \
/180 ; x \oplus y ; z \
/102 ; \'R' x ; \'Q' y ; \'Q' z \
/106 ; \'R' x \cdot \'Q' y ; \'R' y ; \'Q' z \
/72 ; \'R' x \cdot \'Q' y \cdot \'Q' z ; \break \'R' y \cdot \'Q' z ; \'R' z \
\"

\" \t//183. x \oplus ( y \oplus z ) = ( x \oplus y ) \oplus z  \hfill\star \"

\sam183.
\"\p/183.
/\'H' : x : y : z \
/181 ; x ; y ; z \
/182 ; x ; y ; z \
/166 ; x \oplus ( y \oplus z ) \
/166 ; ( x \oplus y ) \oplus z \
/4 ; ( x \oplus y ) \oplus z ; x \oplus ( y \oplus z ) \
\"

\" \t//184. x_1 \le y_1 \and x_2 \le y_2 \imp x_1 + x_2 \le y_1 + y_2 \"

\sam184.
\"\p/184.
/\'H' : x_1 : y_1 : x_2 : y_2 \
/167 ; x_1 ; y_1 : w_1 \
/167 ; x_2 ; y_2 : w_2 \
/72 ; x_1 ; w_1 ; x_2 + w_2 \
/72 ; w_1 ; x_2 ; w_2 \
/98 ; w_1 ; x_2 \
/72 ; x_2 ; \break w_1 ; w_2 \
/72 ; x_1 ; x_2 ; w_1 + w_2 \
/71 ; x_1 + x_2 ; w_1 + w_2 \
\"

\smallskip

With our identification of strings with numbers, 0 is the empty string, 1 is the zero bit, and 2 is
the one bit. For greater readability, introduce new notation for these objects emphasizing their role
as strings.

\" \de//185. \ep = 0 \hfill\star \"

\" \de//186. \u0 = 1 \hfill\star \"

\" \de//187. \u1 = 2 \hfill\star \"

\" \t//188. \u0 \ne \ep \and \u1 \ne \ep \and \u0 \ne \u1 \hfill\star \"

\sam188.
\"\p/188.
/\'H' \
/116 \
/115 \
/185 \
/186 \
/187 \
/3 ; 0 \
/3 ; 1 \
/4 ; 0 ; 1 \
\"

\" \t//189. \'Q' \ep = 1 \and \'R' \ep = 0 \"

\sam189.
\"\p/189.
/\'H' \
/116 \
/115 \
/185 \
/12 ; 1 \
/125 ; 0 \
/130 \
/171 ; 0 ; 1 ; 0 \
\"

\" \t//190. 2 \' is a power of two' \"

\sam190.
\"\p/190.
/\'H' \
/116 \
/115 \
/131 ; 1 \
/126 ; 2 ; 0 \
/117 ; 2 \
/105 ; 2 ; 1 \
\"

\" \t//191. \'Q' \u0 = 2 \and \'R' \u0 = 0 \"

\sam191.
\"\p/191.
/\'H' \
/116 \
/115 \
/186 \
/12 ; 2 \
/190 \
/125 ; 0 \
/125 ; 1 \
/56\fw ; 0 ; 1 \
/122 ; 0 ; 1 ; 2 \
/171 ; 1 ; 2 ; 0 \
\"

\" \t//192. \'Q' \u1 = 2 \and \'R' \u1 = 1 \"

\sam192.
\"\p/192.
/\'H' \
/116 \
/115 \
/187 \
/12 ; 2 ; 0 \
/190 \
/125 ; 1 \
/171 ; 2 ; 2 ; 1 \
\"

\" \t//193. \'Q' x = \'Q' y \and \'R' x = \'R' y \imp x = y \"

\sam193.
\"\p/193.
/\'H' : x : y \
/166 ; x \
/166 ; y \
/4 ; x ; y \
\"

\" \t//194. \ep \oplus x = x \hfill\star \"

\sam194.
\"\p/194.
/\'H' : x \
/193 ; \ep \oplus x ; x \
/189 \
/180 ; \ep ; x \
/117 ; \'Q' x \
/103 ; \'Q' x \
/97 ; \'R' x \
\"

\" \t//195. x \cdot 1 = x \"

\sam195.
\"\p/195.
/\'H' : x \
/117 ; x \
/105 ; x ; 1 \
\"

\" \t//196. x \oplus \ep = x \hfill\star \"

\sam196.
\"\p/196.
/\'H' : x \
/193 ; x \oplus \ep ; x \
/189 \
/180 ; x ; \ep \
/195 ; \'Q' x \
/195 ; \'R' x \
/12 ; \'R' x \
\"

\" \t//197. y \oplus x = z \oplus x \imp y = z \hfill\star \"

\sam197.
\"\p/197.
/\'H' : y : x : z \
/180 ; y ; x \
/180 ; z ; x \
/178 ; x \
/141 ; \'Q' y ; \'Q' x ; \'Q' z \
/119 ; \'R' y \cdot \'Q' x ; \'R' x ; \'R' z \cdot \'Q' x \
/141 ; \'R' y  ; \'Q' x ; \'R' z \
/193 ; y ; z \
\"

\" \t//198. x \cdot y = x \cdot z \and x \ne 0 \imp y = z \"

\sam198.
\"\p/198.
/\'H' : x : y : z \
/105 ; x ; y \
/105 ; x ; z \
/141 ; y ; x ; z \
\"

\" \t//199. x \oplus y = x \oplus z \imp y = z \hfill\star \"

\sam199.
\"\p/199.
/\'H' : x : y : z \
/180 ; x ; y \
/180 ; x ; z \
/178 ; x \
/198 ; \'Q' x ; \'Q' y ; \'Q' z \
/120 ; \'R' x \cdot \'Q' y ; \'R' y ; \'R' z \
/193 ; y ; z \
\"

\" \t//200. x \le 1 \imp x = 0 \or x = 1 \"

\sam200.
\"\p/200.
/\'H' : x \
/116 \
/115 \
/167 ; x ; 1 : w \
/22 ; w \
/12 ; x ; \'P' w \
/4 ; x + \'P' w ; 0 \
/15 ; x ; \'P' w \
\"

\" \t//201. x \ne 0 \imp y \le x \cdot y \"

\sam201.
\"\p/201.
/\'H' : x : y \
/22 ; x \
/104 ; \'P' x ; y \
/98 ; \'P' x \cdot y ; y \
/71 ; y ; \'P' x \cdot y \
\"

\" \t//202. x \cdot y = 1 \imp x = 1 \and y = 1 \"

\sam202.
\"\p/202.
/\'H' : x : y \
/116 \
/115 \
/103 ; y \
/105 ; x ; y \
/103 ; x \
/201 ; x ; y \
/201 ; y ; x \
/200 ; x \
/200 ; y \
/3 ; 0 \
\"

\" \t//203. \'Q' x = 1 \imp x = \ep \"

\sam203.
\"\p/203.
/\'H' : x \
/116 \
/115 \
/185 \
/12 ; 1 \
/166 ; x \
/200 ; \'R' x \
/56\fw ; \'R' x ; 1 \
/4 ; x ; 0 \
\"

\" \t//204. x \oplus y = \ep \imp x = \ep \and y = \ep \hfill\star \"

\sam204.
\"\p/204.
/\'H' : x : y \
/180 ; x ; y \
/189 \
/202 ; \'Q' x ; \'Q' y \
/203 ; x \
/203 ; y \
\"


\" \ar//205. \'Parity' \, 0 = 0 \and \'Parity' \, \'S' x = \'C' ( \'Parity' \, x , 1 , 0 ) \"

\smallskip

$\'Parity' \; x = 0$ expresses that $x$ is even, and $\'Parity' \; x = \'S'0$ that $x$ is odd.

\" \t//206. \'Parity' \, x = 0 \imp \'Parity' \, \'S' x = 1 \"

\sam206.
\"\p/206.
/\'H' : x \
/205 ; x \
/33 ; 1 ; 0 ; 0 \
\"

\" \t//207. \'Parity' \, x = 1 \imp \'Parity' \, \'S' x = 0 \"

\sam207.
\"\p/207.
/\'H' : x \
/116 \
/115 \
/205 ; x \
/33 ; 1 ; 0 ; 0 \
\"

\" \t//208. \'Parity' \, 0 = 0 \and \'Parity' \, 1 = 1 \and \'Parity' \, 2 = 0 \and
\'Parity' \, \ep = 0 \and \'Parity' \, \u0 = 1 \and \'Parity' \, \u1 = 0 \"

\sam208.
\"\p/208.
/\'H' \
/116 \
/115 \
/205 ; 0 \
/205 ; 1 \
/185 \
/186 \
/187 \
/33 ; 1 ; 0; 0 \
\"

\" \t//209. \'Parity' \, x = 0 \or \'Parity' \, x = 1 \"

\sam209.
\"\p/209.
/\'H' : x \
/116 \
/115 \
/205 ; x \
/205 ; \'P' x \
/22 ; x \
/39 ; \'Parity' \, \'P' x \
\"

\" \t//210b. \'Parity' \, x = 0 \imp \'Parity' ( x + 0 ) = \'Parity' \, 0 \"

\sam210b.
\"\p/210b.
/\'H' : x \
/208 \
/12 ; x \
\"

\" \t/210i. [ \'Parity' \, x = 0 \imp \'Parity' ( x + y ) = \'Parity' \, y ] \imp
[ \'Parity' \, x = 0 \imp \'Parity' ( x + \'S' y ) = \'Parity' \, \'S' y ] \"

\sam210i.
\"\p/210i.
/\'H' : x : y \
/12 ; x ; y \
/205 ; x + y \
/205 ; y \
\"

\" \at/210. \'Parity' \, x = 0 \imp \'Parity' ( x + y ) = \'Parity' \, y \"

\" \t//211b. \'Parity' \, x = 0 \imp \'Parity' ( x \cdot 0 ) = 0 \"

\sam211b.
\"\p/211b.
/\'H' : x \
/100 ; x \
/205 \
\"

\"\t/211i. [ \'Parity' \, x = 0 \imp \'Parity' ( x \cdot y ) = 0 ] \imp
[ \'Parity' \, x = 0 \imp \'Parity' ( x \cdot \'S' y ) = 0 ] \"

\sam211i.
\"\p/211i.
/\'H' : x : y \
/99 ; x ; y \
/210 ; x ; x \cdot y \
\"

\" \at/211. \'Parity' \, x = 0 \imp \'Parity' ( x \cdot y ) = 0 \"

\" \t//212. \'Parity' ( x \cdot 2 + y ) = \'Parity' \, y \"

\sam212.
\"\p/212.
/H : x : y \
/208 \
/105 ; x ; 2 \
/211 ; 2 ; x \
/210 ; 2 \cdot x ; y \
\"

\" \t//213. \'Parity' ( x \cdot 2 ) = 0 \"

\sam213.
\"\p/213.
/\'H' : x \
/212 ; x ; 0 \
/208 \
/12 ; x \cdot 2 \
\"

\" \t//214. \'Parity' ( x \cdot 2 + 1 ) = 1 \"

\sam214.
\"\p/214.
/\'H' : x \
/212 ; x ; 1 \
/208 \
\"

\" \t//215. \'Parity' \, \'R' ( x \oplus \u0 ) = 0 \"

\sam215.
\"\p/215.
/\'H' : x \
/180 ; x ; \u0 \
/191 \
/212 ; \'R' x ; 0 \
/205 \
\"

\" \t//216. \'Parity' \, \'R' ( x \oplus \u1 ) = 1 \"

\sam216.
\"\p/216.
/\'H' : x \
/180 ; x ; \u1 \
/192 \
/212 ; \'R' x ; 1 \
/208 \
\"

\" \t//217. x \oplus \u0 \ne y \oplus \u1 \hfill\star \"

\sam217.
\"\p/217.
/\'H' : x : y \
/116 \
/115 \
/215 ; x \
/216 ; y \
/3 ; 0 \
\"

\smallskip

{\it \#12. Construct the\/ unary {\rm PR} function symbol\/ {\rm Chop} that deletes the last bit, if any,
of a string.}


\" \ar//218. \'Half' \; 0 = 0 \and \'Half' \; \'S' x =
\'C' ( \'Parity' \; x , \'Half' \; x , \'S' \, \'Half' \; x ) \"

\smallskip

$\'Half' \; x $ expresses $\lfloor x\slash{2}\rfloor$. It deletes the rightmost binary bit of $x$
(if $x$ is not 0 or 1).

\" \t//219. \'Parity' \, x = 0 \imp \'Half' \; \'S' x = \'Half' \; x \"

\sam219.
\"\p/219.
/\'H' : x \
/218 ; x \
/33 ; \'Half' \, x ; \'S' \, \'Half' x \
\"

\" \t//220. \'Parity' \, x = 1 \imp \'Half' \; \'S' x = \'S' \, \'Half' \; x \"

\sam220.
\"\p/220.
/\'H' : x \
/116 \
/115 \
/218 ; x \
/33 ; \'Half' \, x ; \'S' \, \'Half' x; 0 \
\"

\" \t//221. \'Half' \; 0 = 0 \and \'Half' \; 1 = 0 \and \'Half' \; 2 = 1 \and
\'Half' \; \'S' 2 = 1 \"

\sam221.
\"\p/221.
/\'H' \
/116 \
/115 \
/218 ; 0 \
/218 ; 1 \
/218 ; 2 \
/208 \
/219 ; x \
/33 ; 0 ; 1 ; 0 \
/33 ; 0 ; 2 ; 1 \
/33 ; 1 ; 2 ; 0 \
\"

\" \d//222. \'p'_{/222} ( x ) \iff [ \'Parity' \, x = 0 \imp x = 2 \cdot \'Half' \; x ] \and
[ \'Parity' \, x = 1 \imp x = 2 \cdot \'Half' \; x + 1 ] \"

\" \t//223b. \'p'_{/222} ( 0 ) \"

\sam223b.
\"\p/223b.
/\'H' \
/116 \
/115 \
/222\bw ; 0 \
/205 \
/221 \
/99 ; 2 \
/3 ; 0 \
\"

\" \t/223i. \'p'_{/222} ( x ) \imp \'p'_{/222} ( \'S' x ) \"

\sam223i.
\"\p/223i.
/\'H' : x \
/116 \
/115 \
/222\fw ; x \
/222\bw ; \'S' x \
/? \'Parity' \, x = 1 \
/209 ; x \
/3 ; 0 \
/12 ; 2 \cdot \'Half' \, x ; 0 \
/206 ; x \
/207 ; x \
/219 ; x \
/220 ; x \
/12 ; 2 \cdot \'Half' \, x ; 1 \
/100 ; 2 ; \'Half' \, x \
\"

\" \at/223. \'p'_{/222} ( x ) \"

\" \t//224. \'Parity' \, x = 0 \imp x = 2 \cdot \'Half' \; x \"

\sam224.
\"\p/224.
/\'H' : x \
/223 ; x \
/222\fw ; x \
\"

\" \t//225. \'Parity' \, x = 1 \imp x = 2 \cdot \'Half' \; x + 1 \"

\sam225.
\"\p/225.
/\'H' : x \
/223 ; x \
/222\fw ; x \
\"

\" \t//226. x \le 2 \cdot \'Half' \; x + 1 \"

\sam226.
\"\p/226.
/\'H' : x \
/224 ; x \
/225 ; x \
/71 ; 2 \cdot \'Half' \, x ; 1 \
/209 ; x \
/60 ; x \
\"

\" \t//227. \'Half' \; x = 0 \imp x = 0 \or x = 1 \"

\sam227.
\"\p/227.
/\'H' : x \
/226 ; x \
/100 ; 2 \
/97 ; 1 \
/200 ; x \
\"

\" \t//228. x \ne \ep \imp \'Half' \; \'Q' x + \'Half' \; \'R' x \ne 0 \"

\sam228.
\"\p/228.
/\'H' : x \
/15 ; \'Half' \, \'Q' x ; \'Half' \, \'R' x \
/227 ; \'Q' x \
/178 ; x \
/203 ; x
\"


\noin e/229. $\ \'Chop' \, x = \'P' ( \'Half' \; \'Q' x + \'Half' \; \'R' x ) $

\smallskip

$\'Chop' \, x$ deletes the rightmost bit, if any, of the string $x$. It is the string analogue of P.

\" \t//230. \'Chop' \, \ep = \ep \hfill\star \"

\sam230.
\"\p/230.
/\'H' \
/229 ; \ep \
/189 \
/185 \
/221 \
/12 ; 0 \
/16 \
\"

\" \t//231. x \ne \ep \imp \'S' \, \'Chop' \, x = \'Half' \; \'Q' x + \'Half' \; \'R' x \"

\sam231.
\"\p/231.
/\'H' : x \
/229 ; x \
/228 ; x \
/22 ; \'Half' \, \'Q' x + \'Half' \, \'R' x \
\"

\" \t//232. 2 \cdot \'Half' \; x \le x \"

\sam232.
\"\p/232.
/\'H' : x \
/224 ; x \
/225 ; x \
/60 ; x \
/71 ; 2 \cdot \'Half' \, x ; 1 \
/209 ; x \
\"

\" \t//233. \'Parity' \, q = 0 \and r \lt q \imp \'Half' \; r \lt \'Half' \; q \"

\sam233.
\"\p/233.
/\'H' : q : r \
/79 ; \'Half' \, q ; \'Half' \, r \
/224 ; q \
/118 ; \'Half' \, q \
/184 ; \'Half' \, q ; \'Half' \, r ; \'Half' \, q ; \'Half' \, r \
/232 ; r \
/118 ; \'Half' \, r \
/73 ; q ; \'Half' \, r + \'Half' \, r ; r \
/80 ; r ; q \
\"

\" \t//234. \'Parity' ( 2 \up \'S' x ) = 0 \"

\sam234.
\"\p/234.
/\'H' : x \
/126 ; 2 ; x \
/208 \
/211 ; 2 ; 2 \up x \
\"

\" \t//235. \'Half' ( 2 \up \'S' x ) = 2 \up x \"

\sam235.
\"\p/235.
/\'H' : x \
/116 \
/115 \
/234 ; x \
/224 ; 2 \up \'S' x \
/126 ; 2 ; x \
/3 ; 1 \
/198 ; 2 ; \'Half' ( 2 \up \'S' x ) ; 2 \up x \
\"

\" \t//236. q \' is a power of two' \and q \ne 1 \imp \'Parity' \; q = 0 \"

\sam236.
\" \p/236.
/\'H' : q \
/129\fw ; q : x \
/126 ; 2 ; x \
/22 ; x \
/134 \
/234 ; \'P' x \
\"

\" \t//237. x \ne 0 \imp \'Half' ( 2 \up x ) \' is a power of two' \"

\sam237.
\"\p/237.
/\'H' ; x \
/116 \
/115 \
/22 ; x \
/235 ; \'P' x \
/131 ; \'P' x \
\"

\" \t//238. x \ne \ep \imp \'Half' \; \'Q' x \' is a power of two' \"

\sam238.
\"\p/238.
/\'H' : x \
/158 ; x \
/116 \
/115 \
/129\fw ; \'Q' x : y \
/126 ; 2 \
/203 ; x \
/237 ; y \
\"

\" \t//239. x \ne \ep \imp \'Q' \, \'Chop' \, x = \'Half' \; \'Q' x \and \'R' \, \'Chop' \, x = \'Half' \; \'R' x \"

\sam239.
\"\p/239.
/\'H' : x \
/171 ; \'Chop' \, x ; \'Half' \, \'Q' x ; \'Half' \, \'R' x \
/231 ; x \
/238 ; x \
/236 ; \'Q' x \
/166 ; x \
/203 ; x \
/233 ; \'Q' \, x ; \'R' x \
\"

\" \t//240. \'Chop' \, \u0 = \ep \and \'Chop' \, \u1 = \ep \hfill\star \"

\sam240.
\"\p/240.
/\'H' \
/239 ; \u0 \
/239 ; \u1 \
/188 \
/191 \
/192 \
/221 \
/203 ; \'Chop' \, \u0 \
/203 ; \'Chop' \, \u1 \
\"

\" \t//241b. \'Half' ( 0 \cdot 2 ) = 0 \and \'Half' ( 0 \cdot 2 + 1 ) = 0 \"

\sam241b.
\"\p/241b.
/\'H' \
/103 ; 2 \
/97 ; 1 \
/221 \
\"

\" \t/241i. \'Half' ( x \cdot 2 + 1 ) = x \and \'Half' ( x \cdot 2 ) = x \imp
\'Half' ( \'S' x \cdot 2 + 1 ) = \'S' x \and \'Half' ( \'S' x \cdot 2 ) = \'S' x \"

\sam241i.
\"\p/241i.
/\'H' : x \
/116 \
/115 \
/104 ; x ; 2 \
/12 ; x \cdot 2 ; 1 \
/212 ; x ; 1 \
/208 \
/220 ; x \cdot 2 + 1 \
/12 ; x \cdot 2 + 2 ; 0 \
/12 ; \'S' ( x \cdot 2 + 2 ) \
/212 ; x ; 2 \
/219 ; x \cdot 2 + 2 \
\"

\" \at/241. \'Half' ( x \cdot 2 + 1 ) = x \and \'Half' ( x \cdot 2 ) = x \"

\" \t//242. x \oplus \u0 \ne \ep \and x \oplus \u1 \ne \ep \"

\sam242.
\"\p/242.
/\'H' : x \
/204 ; x ; \u0 \
/204 ; x ; \u1 \
/188 \
\"

\" \t//243. \'Chop' ( x \oplus \u0 ) = x \hfill\star \"

\sam243.
\"\p/243.
/\'H' : x \
/242 ; x \
/239 ; x \oplus \u0 \
/180 ; x ; \u0 \
/191 \
/241 ; \'Q' x \
/241 ; \'R' x \
/12 ; \'R' x \cdot 2 \
/193 ; x ; \'Chop' ( x \oplus \u0 ) \
\"

\" \t//244. \'Chop' ( x \oplus \u1 ) = x \hfill\star \"

\sam244.
\"\p/244.
/\'H' : x \
/242 ; x \
/239 ; x \oplus \u1 \
/180 ; x ; \u1 \
/192 \
/241 ; \'Q' x \
/241 ; \'R' x \
/12 ; \'R' x \cdot 2 \
/193 ; x ; \'Chop' ( x \oplus \u1 ) \
\"

\" \t//245. x \ne \ep \imp \'Q' ( \'Chop' \, x \oplus \u0 ) = \'Q' x
\and \'Q' ( \'Chop' \, x \oplus \u1 ) = \'Q' x \"

\sam245.
\"\p/245.
/\'H' : x \
/180 ; \'Chop' \, x ; \u0 \
/180 ; \'Chop' \, x ; \u1 \
/191 \
/192 \
/158 ; x \
/203 ; x \
/236 ; \'Q' x \
/224 ; \'Q' x \
/239 ; x \
/105 ; 2 ; \'Q' \, \'Chop' \, x \
\"

\" \t//246. \'Parity' \, x = \'Parity' \, y \and \'Half' \; x = \'Half' \; y \imp x = y \"

\sam246.
\"\p/246.
/\'H' : x : y \
/224 ; x \
/224 ; y \
/225 ; x \
/225 ; y \
/221 ; x \
/209 ; x \
/? \'Parity' \, x = 0 \
\"

\" \t//247. x \ne \ep \imp \'R' ( \'Chop' x \oplus \u0 ) = \'Half' \, \'R' \, x \cdot 2 \"

\sam247.
\"\p/247.
/\'H' : x \
/180 ; \'Chop' \, x ; \u0 \
/191 \
/12 ; \'R' \, \'Chop' \, x \cdot 2 \
/239 ; x \
\"

\" \t//248. x \ne \ep \and \'Parity' \, \'R' x = 0 \imp \'R' ( \'Chop' \, x \oplus \u0 ) = \'R' x \"

\sam248.
\"\p/248.
/\'H' : x \
/247 ; x \
/241 ; \'Half' \, \'R' x \
/213 ; \'Half' \, \'R' x \
/246 ; \'R' ( \'Chop' \, x \oplus \u0 ) ; \'R' x \
\"

\" \t//249. x \ne \ep \imp \'R' ( \'Chop' x \oplus \u1 ) = \'Half' \; \'R' x \cdot 2 + 1 \"

\sam249.
\"\p/249.
/\'H' : x \
/180 ; \'Chop' \, x ; \u1 \
/192 \
/239 ; x \
\"

\" \t//250. x \ne \ep \and \'Parity' \, \'R' x = 1 \imp \'R' ( \'Chop' \, x \oplus \u1 ) = \'R' x \"

\sam250.
\"\p/250.
/\'H' : x \
/249 ; x \
/241 ; \'Half' \, \'R' x \
/214 ; \'Half' \, \'R' x \
/246 ; \'R' ( \'Chop' \, x \oplus \u1 ) ; \'R' x \
\"

\" \t//251. x \ne \ep \and \'Parity' \; \'R' x = 0 \imp x = \'Chop' \, x \oplus \u0 \"

\sam251.
\"\p/251.
/\'H' : x \
/245 ; x \
/248 ; x \
/193 ; x ; \'Chop' \, x \oplus \u0 \
\"

\" \t//252. x \ne \ep \and \'Parity' \; \'R' x = 1 \imp x = \'Chop' \, x \oplus \u1 \"

\sam252.
\"\p/252.
/\'H' : x \
/245 ; x \
/250 ; x \
/193 ; x ; \'Chop' \, x \oplus \u1 \
\"

\" \t//253. x \ne \ep \imp x = \'Chop' \, x \oplus \u0 \or x = \'Chop' \, x \oplus \u1 \hfill{\star}{\star} \"

\sam253.
\"\p/253.
/\'H' : x \
/251 ; x \
/252 ; x \
/209 ; \'R' x \
\"

\" \t//254. y \ne \ep \imp \'Chop' ( x \oplus y ) = x \oplus \'Chop' \, y \hfill\star \"

\sam254.
\"\p/254.
/\'H' : y : x \
/253 ; y \
/183 ; x ; \'Chop' \, y ; \u0 \
/183 ; x ; \'Chop' \, y ; \u1 \
/243 ; x \oplus \'Chop' \, y \
/244 ; x \oplus \'Chop' \, y \
\"

\smallbreak
\goodbreak

{\it \#13. Establish string recursion and string induction.}

\everymath = {\rm} \everydisplay = {\rm}

\Pr//  8. {\it Let\/ $\vec x$, \ro y, and\/ \ro z be distinct, let\/ \ro a contain no variables other
than those in\/ $\vec x$, and let\/ \ro b and\/ \ro c contain no variables other than those in\/ $\vec x$,
\ro y, and\/ \ro z. Define\/ \ro f by primitive recursion by cases:
$$\eqalign{&f(\vec x,0)=a\and f(\vec x,Sy)= {}\cr
&C\big(\chi[Sy=Chop\,Sy\oplus\u0],b_z(Chop\,Sy),
C(\chi[Sy=Chop\,Sy\oplus\u1],c_z(Chop\,Sy),0)\big)\cr}\leqno(// 6)$$
Then
$$f(\vec x,\ep)=a\and f(\vec x, w\oplus\u0)=b_z(w)\and f(\vec x,w\oplus\u1)=c_z(w)\leqno(// 7)$$}

\vskip-15pt

\pf We have $f(\vec x,\ep)=a$ by \Rf{/185}. By Proposition/  4 of $\S/   4$,
$$\leqalignno{
Sy=Chop\,Sy\oplus\u0 &\imp f(\vec x,Sy)=b_z(Chop\,Sy)&(// 8)\cr
\Neg[Sy=Chop\,Sy\oplus\u0]\and[Sy=Chop\,Sy\oplus\u1] &\imp f(\vec x,Sy)=c_z(Chop\,Sy)&(// 9)\cr
\Neg[Sy=Chop\,Sy\oplus\u0]\and\Neg[Sy=Chop\,Sy\oplus\u1] &\imp f(\vec x,Sy)=0&(// 10)\cr
}$$
By $\rf/217;Chop\,Sy;Chop\,Sy$ and $(/ 9)$,
$$Sy=Chop\,Sy\oplus\u1 \imp f(\vec x,Sy)=c_z(Chop\,Sy)\leqno(/ 11)$$
Note that $(/ 10)$, the ``else'' clause, is irrelevant, since its hypothesis cannot hold,
by \rf/253;Sy together with \rf/3;y and \Rf{/185}.
Now consider $(/ 8);P(w\oplus\u0)$:
$$SP(w\oplus\u0)=Chop\,SP(w\oplus\u0)\oplus\u0\imp f(\vec x,SP(w\oplus\u0))=b_z(Chop\,SP(w\oplus\u0))$$
We have $w\oplus\u0\ne0$ (by \Rf{/185}\ and $\rf/204;w;\u0$ and \Rf{/188}), so $SP(w\oplus\u0)=w\oplus\u0$ by
$\rf/22;w\oplus\u0$. Also, $Chop(w\oplus\u0)=w$ by \rf/243;w. Therefore
$$w\oplus\u0=w\oplus\u0\imp f(\vec x,w\oplus\u0)=b_z(w)$$
and so $f(\vec x,w\oplus\u0)=b_z(w)$.
The derivation of $f(\vec x,w\oplus\u1)=c_z(w)$ from $(/ 11);SP(w\oplus\u1)$ is entirely similar.
Hence $(/ 7)$. \bul

\smallskip

A \ii string recursion.\/ is a primitive recursion of the form $(/ 6)$, but string recursions will be
introduced simply by $(/ 7)$.

\everymath={} \everydisplay={}

\" \t//255. x \lt x \oplus \u0 \"

\sam255.
\"\p/255.
/\'H' : x \
/172 ; x ; \u0 \
/191 \
/12 ; \'S' x \cdot 2 \
/105 ; \'S' x ; 2 \
/118 ; \'S' x \
/12 ; \'S' x ; x \
/16 ; \'S' x + x \
/71 ; \'S' x ; x \
/125 ; x \
/144 ; x ; \'S' x ; \'S' x + x \
\"

\" \t//256. x + 1 = \'S' x \"

\sam256.
\"\p/256.
/\'H' : x \
/115 \
/12 ; x ; 0 \
\"

\" \t//257. x \lt x \oplus \u1 \"

\sam257.
\"\p/257.
/\'H' : x \
/172 ; x ; \u1 \
/192 \
/256 ; \'S' x \cdot 2 \
/16 ; \'S' x \cdot 2 \
/105 ; \'S' x ; 2 \
/118 ; \'S' x \
/125 ; x \
/71 ; \'S' x ; \'S' x \
/144 ; x ; \'S' x ; \'S' x + \'S' x \
\"

\" \t//258. x \ne \ep \imp \'Chop' \, x \lt x \"

\sam258.
\"\p/258.
/\'H' : x \
/253 ; x \
/255 ; \'Chop' \, x \
/257 ; \'Chop' \, x \
\"

\goodbreak

\everymath={\rm} \everydisplay={\rm}

\Pr//  9. {\it If\ $\vdash A_x(\ep)$ and\/ $\vdash A_x(x')\impP A_x(x'\oplus\u0)$ and\/
$\vdash A_x(x')\impP A_x(x'\oplus\u1)$, then\/ $\vdash A_x(x)$.}

\pf Suppose that
$$\leqalignno{
&A_x(\ep)&(// 12)\cr
&A_x(x')\imp A_x(x'\oplus\u0)&(// 13)\cr
&A_x(x')\imp A_x(x'\oplus\u1)&(// 14)\cr
}$$
and use the least number principle (Proposition/  6 of $\S/   5$).
Suppose that there is a least counterexample~z to~A:
$$\Neg A_x(z) \and \forall y\ltT z [A_x(y)] \leqno(// 15)$$
We have $z\ne\ep$ by $(/ 12)$. Then, by \rf/253;z,
$$z=Chop\,z\oplus\u0\or z=Chop\,z\oplus\u1\leqno(// 16)$$
By $(/ 15);z;Chop\,z$ and \rf/258;z,
$$A_x(Chop\,z)$$
The first alternative of $(/ 16)$ does not hold, by $(/ 13);Chop\,z$, and the second does not
hold, by $(/ 14);Chop\,z$. Hence $A_x(x)$ by the least number principle. \bul

\smallskip

The derived rule of inference, from $A_x(\ep)$ and $A_x(x')\impP A_x(x'\oplus\u0)$ and
$A_x(x')\impP A_x(x'\oplus\u1)$ infer $A_x(x)$, is \ii string induction..

\everymath={} \everydisplay={}

\smallskip

{\it \#14. Construct\/ {\rm Length} by string recursion, and introduce $\preceq$ (shorter than) and $\prec$
(strictly shorter than). Prove $x\preceq y\impP x\le 2\cdot\ro Qx$ and $x\prec y\impP x\lt y$. The first will
enable bounded quantifiers with bounding symbol $\preceq$ rather than $\le$, and the second will enable the
use of the least number principle for strings (``the shortest string principle''). Express the i'th bit of
a string.}


\" \ar//259. \'Length' \, \ep = 0 \and \'Length' ( x \oplus \u0 ) = \'S' \, \'Length' \, x
\and \'Length' ( x \oplus \u1 ) = \'S' \, \'Length' \, x \hfill\star \"

\" \t//260b. \'Length' ( x \oplus \ep ) = \'Length' \, x + \'Length' \, \ep \"

\sam260b.
\"\p/260b.
/\'H' : x \
/196 ; x \
/259 \
/12 ; \'Length' \, x \
\"

\" \t/260ij. \'Length' ( x \oplus y ) = \'Length' \, x + \'Length' \, y \imp
\'Length' \big ( x \oplus ( y \oplus \u0 ) \big ) = \'Length' \, x + \'Length' ( y \oplus \u0 ) \and
\'Length' \big ( x \oplus ( y \oplus \u1 ) \big ) = \'Length' \, x + \'Length' ( y \oplus \u1 ) \"

\sam260ij.
\"\p/260ij.
/\'H' : x : y \
/259 ; y \
/183 ; x ; y ; \u0 \
/183 ; x ; y ; \u1 \
/259 ; x \oplus y \
/12 ; \'Length' \, x ; \'Length' \, y \
\"

\" \at/260. \'Length' ( x \oplus  y ) = \'Length' \, x + \'Length' \, y \"

\smallskip

When t$\xi$ immediately follows t$\xi b$, t$\xi i$, and t$\xi j$, it is an inference by string induction.
Sometimes, as here, the two string induction steps are combined into a single formula~t$\xi ij$.

\" \t//261b. \'Q' \ep = 2 \up ( \'Length' \, \ep ) \"

\sam261b.
\"\p/261b.
/\'H' \
/259 \
/126 ; 2 \
/189 \
\"

\" \t/261ij. \'Q' x = 2 \up \'Length' \, x \imp
\'Q' ( x \oplus \u0 ) = 2 \up \'Length' ( x \oplus \u0 ) \and
\'Q' ( x \oplus \u1 ) = 2 \up \'Length' ( x \oplus \u1 ) \"

\sam261ij.
\"\p/261ij.
/\'H' : x \
/259 ; x \
/126 ; 2 ; \'Length' \, x \
/180 ; x ; \u0 \
/180 ; x ; \u1 \
/191 \
/192 \
/105 ; \'Q' x ; 2 \
\"

\" \at/261. \'Q' x = 2 \up \'Length' \, x \"

\" \t//262. x \le \'S' y \imp x \le y \or x = \'S' y \"

\sam262.
\"\p/262.
/\'H' : x : y \
/68 ; x ; \'S' y \
/22 ; \'S' y \minus x \
/12 ; x ; \'P' ( \'S' y \minus x ) \
/4 ; y ; x + \'P' ( \'S' y \minus x ) \
/71 ; x ; \'P' ( \'S' y \minus x ) \
\"

\" \t//263b. x \le 0 \imp 2 \up x \le 2 \up 0 \"

\sam263b.
\"\p/263b.
/\'H' : x \
/57 ; x \
/60 ; 2 \up 0 \
\"

\" \t/263i. [ x \le y \imp 2 \up x \le 2 \up y ] \imp [ x \le \'S' y \imp 2 \up x \le 2 \up \'S' y ] \"

\sam263i.
\"\p/263i.
/\'H' : x : y \
/262 ; x ; y \
/60 ; 2 \up x \
/126 ; 2 ; y \
/118 ; 2 \up y \
/71 ; 2 \up y ; 2 \up y \
/73 ; 2 \up x ; 2 \up y ; 2 \cdot ( 2 \up y ) \
\"

\" \at/263. x \le y \imp 2 \up x \le 2 \up y \"

\" \d//264. x \preceq y \iff \'Length' \, x \le \'Length' \, y \"

\" \d//265. x \prec y \iff \'Length' \, x \lt \'Length' \, y \"

\" \t//266. x \le y \imp 2 \cdot x \le 2 \cdot y \"

\sam266.
\"\p/266.
/\'H' : x : y \
/118 ; x \
/118 ; y \
/184 ; x ; y ; x ; y \
\"

\" \t//267. x \preceq y \imp x \le 2 \cdot \'Q' y \hfill\star \"

\sam267.
\"\p/267.
/\'H' : x : y \
/264\fw ; x ; y \
/263 ; \'Length' \, x ; \'Length' \, y \
/261 ; x \
/261 ; y \
/158 ; x \
/63 ; x \
/56\fw ; \'S' x ; 2 \cdot \'Q' x \
/73 ; x ; \'S' x ; 2 \cdot \'Q' x \
/266 ; \'Q' x ; \'Q' y \
/73 ; x ; 2 \cdot \'Q' x ; 2 \cdot \'Q' y \
\"

\smallskip

Let s be a binary predicate symbol, written in infix notation. If for some term d containing no variable other
than~$y$ we have
$\vdash x \; \ro s \; y \impP x \le \ro d$, then we call s a \ii bounding symbol.\/ with \ii bound.\/ d.
Thus $\preceq$ is a bounding symbol with bound $2\cdot\'Q' y$ by t/267. For a bounding symbol~s with
bound~d let $\rm\exists x\,s\,b\, A$ abbreviate
$\rm\exists x\leE d_{\mit y}(b)[x \;s\;b \andD A]$ and $\rm\forall x\,s\,b\, A$ abbreviate
$\rm\forall x\leE d_{\mit y}(b)[x \;s\;b \impP A]$. (The bound d is not unique, but we choose one, and it is
implicit in the notations $\rm\exists x\,s\,b\, A$ and $\rm\forall x\,s\,b\, A$.)
We have already used these abbreviations for the bounding symbol~$\lt$ (with bound $y$).

\" \t//268. \'S' \'R' x \lt \'Q' x \imp \'Q' \'S' x = \'Q' x \"

\sam268.
\"\p/268.
/\'H' : x \
/166 ; x \
/12 ; \'Q' x ; \'R' x \
/171 ; \'S' x ; \'Q' x ; \'S' \'R' x \
\"

\" \t//269. q \' is a power of two' \imp 0 \lt q \"

\sam269.
\"\p/269.
/\'H' : q \
/134 \
/58 ; q \
/56\bw ; 0 ; q \
\"

\" \t//270. \'S' \'R' x = \'Q' x \imp \'Q' \'S' x = 2 \cdot \'Q' x \"

\sam270.
\"\p/270.
/\'H' : x \
/166 ; x \
/190 \
/177 ; 2 ; \'Q' x \
/12 ; \'Q' x ; \'R' x \
/118 ; \'Q' x \
/269 ; 2 \cdot \'Q' x \
/12 ; 2 \cdot \'Q' x \
/171 ; \'S' x ; 2 \cdot \'Q' x ; 0 \
\"

\" \t//271. \'Q' x \le \'Q' \'S' x \"

\sam271.
\"\p/271.
/\'H' : x \
/60 ; \'Q' x \
/268 ; x \
/166 ; x \
/270 ; x \
/118 ; \'Q' x \
/71 ; \'Q' x ; \'Q' x \
/114 ; \'R' x ; \'Q' x \
/56\bw ; \'S' \'R' x ; \'Q' x \
\"

\" \t//272b. x \le 0 \imp \'Q' x \le \'Q' 0 \"

\sam272b.
\"\p/272b.
/\'H' : x \
/57 ; x \
/60 ; \'Q' 0 \
\"

\" \t/272i. [ x \le y \imp \'Q' x \le \'Q' y ] \imp
[ x \le \'S' y \imp \'Q' x \le \'Q' \'S' y ] \"

\sam272i.
\"\p/272i.
/\'H' : x : y \
/110 ; x ; y \
/60 ; \'Q' \'S' y \
/271 ; y \
/73 ; \'Q' x ; \'Q' y ; \'Q' \'S' y \
\"

\" \at/272. x \le y \imp \'Q' x \le \'Q' y \"

\" \t//273. x \prec y \imp x \lt y \hfill\star \"

\smallskip

Since $\'Length'\, x \lt \'Length' \,y$, we have $\'Q'x\lt\'Q'y$ by t/261 and t/145.
If $\neg x\lt y$, then $y\le x$,
so $\'Q'y\le\'Q'x$ by t/272, a contradiction. It is essential to have the strict bound~$y$ itself on~$\prec$ to
apply the least number principle to strings.

\sam273.
\"\p/273.
/\'H' : x : y \
/265\fw ; x ; y \
/145 ; \'Length' \, x ; \'Length' \, y \
/261 ; x \
/261 ; y \
/272 ; y ; x \
/80 ; \'Q' x ; \'Q' y \
/79 ; y ; x \
\"

\" \t//274. \ep \preceq x \and x \preceq x \"

\sam274.
\"\p/274.
/\'H' : x \
/264\bw ; \ep ; x \
/264\bw ; x ; x \
/259 \
/58 ; \'Length' \, x \
/60 ; \'Length' \, x \
\"

\" \t//275. x \preceq y \and y \preceq z \imp x \preceq z \"

\sam275.
\"\p/275.
/\'H' : x : y : z \
/264\fw ; x ; y \
/264\fw ; y ; z \
/73 ; \'Length' \, x ; \'Length' \, y ; \'Length' \, z \
/264\bw ; x ; z \
\"

\" \t//276. x \oplus y \preceq y \oplus x \"

\sam276.
\"\p/276.
/\'H' : x : y \
/260 ; x ; y \
/260 ; y ; x \
/98 ; \'Length' \, x ; \'Length' \, y \
/60 ; \'Length' \, x + \'Length' \, y \
/264\bw ; x \oplus y ; y \oplus x \
\"

\" \t//277. x \preceq x \oplus y \and y \preceq x \oplus y \"

\sam277.
\"\p/277.
/\'H' : x : y \
/264\bw ; x ; x \oplus y \
/264\bw ; y ; x \oplus y \
/260 ; x ; y \
/71 ; \'Length' \, x ; \'Length' \, y \
/98 ; \'Length' \, x ; \'Length' \, y \
/71 ; \break \'Length' \, y ; \'Length' \, x \
\"

\" \t//278. \'Length' \, \u0 = 1 \and \'Length' \, \u1 = 1 \"

\sam83.
\"\p/278.
/\'H' \
/259 ; \ep \
/194 ; \u0 \
/194 ; \u1 \
/115 \
\"

\" \t//000b. x_1 \oplus \ep = x_2 \oplus \ep \imp x_1 = x_2 \"

sam000b.
\"\p/000b.
/\'H' : x_1 : x_2 \
/196 ; x_1 \
/196 ; x_2 \
\"

\" \t/000i. [ x_1 \oplus y = x_2 \oplus y \imp x_1 = x_2 ] \imp
[ x_1 \oplus y \oplus \u0 = x_2 \oplus y \oplus \u0 \imp x_1 = x_2 ] \"

\sam000i.
\"\p/000i.
/\'H' : x_1 : y : x_2 \
/183 ; x_1 ; y ; \u0 \
/183 ; x_2 ; y ; \u0 \
/243 ; x_1 \oplus y \
/243 ; x_2 \oplus y \
\"

\" \t/000j. [ x_1 \oplus y = x_2 \oplus y \imp x_1 = x_2 ] \imp
[ x_1 \oplus y \oplus \u1 = x_2 \oplus y \oplus \u1 \imp x_1 = x_2 ] \"

\sam000j.
\"\p/000j.
/\'H' : x_1 : y : x_2 \
/183 ; x_1 ; y ; \u1 \
/183 ; x_2 ; y ; \u1 \
/244 ; x_1 \oplus y \
/244 ; x_2 \oplus y \
\"

\" \at/000. x_1 \oplus y = x_2 \oplus y \imp x_1 = x_2 \"

\" \t//001. \'Length' \, x = 0 \imp x = \ep \"

\sam001.
\"\p/001.
/\'H' : x \
/253 ; x \
/260 ; \'Chop' \, x ; \u0 \
/260 ; \'Chop' \, x ; \u1 \
/278 \
/15 ; \'Length' \, \'Chop' \, x ; 1 \
/115 \
/3 ; 0 \
\"

\" \d//288. x \' ends with ' c \iff \exists a \preceqQ x [ x = a \oplus c ] \"

\" \t//289. a \oplus c \' ends with ' c \"

\sam289.
\"\p/289.
/\'H' : a : c \
/277 ; a ; c \
/288\bw ; a \oplus c ; c ; a \
/5 ; a \oplus c \
\"

\" \d//290. x \' begins with ' a \iff \exists c \preceqQ x [ x = a \oplus c ] \"

\" \t//291. a \oplus c \' begins with ' a \"

\sam291.
\"\p/291.
/\'H' : a : c \
/277 ; a ; c \
/290\bw ; a \oplus c ; a ; c \
/5 ; a \oplus c \
\"

\" \t//292. x \' begins with ' a \imp x \oplus y \' begins with ' a \"

\sam292.
\"\p/292.
/\'H' : x : a : y \
/290\fw ; x ; a : c \
/183 ; a ; c ; y \
/291 ; a ; c \oplus y \
\"

\" \t//293. y \' ends with ' c \imp x \oplus y \' ends with ' c \"

\sam293.
\"\p/293.
/\'H' : y : c : x \
/288\fw ; y ; c : a \
/183 ; x ; a ; c \
/289 ; x \oplus a ; c \
\"

\" \t//294. x \' ends with ' x \and x \' ends with ' \ep \and x \' begins with ' x \and x \' begins with ' \ep \"

\sam294.
\"\p/294.
/\'H' : x \
/194 ; x \
/196 ; x \
/289 ; \ep ; x \
/289 ; x ; \ep \
/291 ; \ep ; x \
/291 ; x ; \ep \
\"

\" \t//002. \ep \' begins with ' x \or \ep \' ends with ' x \imp x = \ep \"

\sam002.
\"\p/002.
/\'H' : x \
/290\fw ; \ep ; x : c \
/288\fw ; \ep ; x : a \
/204 ; x ; c \
/204 ; a ; x \
\"

\" \t//005. x_1 + y_1 = x_2 + y_2 \and x_1 \le x_2 \imp y_2 \le y_1 \"

\sam005.
\"\p/005.
/\'H' : x_1 : y_1 : x_2 : y_2 \
/68 ; x_1 ; x_2 \
/72 ; x_1 ; x_2 \minus x_1 ; y_2 \
/120 ; x_1 ; y_1 ; ( x_2 \minus x_1 ) + y_2 \
/98 ; x_2 \minus x_1 ; y_2 \
/71 ; y_2 ; x_2 \minus x_1 \
\"

\" \t//006. a_1 \oplus c_1 = a_2 \oplus c_2 \and a_1 \preceq a_2 \imp c_2 \preceq c_1 \"

\sam006.
\"\p/006.
/\'H' : a_1 : c_1 : a_2 : c_2 \
/264\fw ; a_1 ; a_2 \
/264\bw ; c_2 ; c_1 \
/260 ; a_1 ; c_1 \
/260 ; a_2 ; c_2 \
/005 ; \'Length' \, a_1 ; \'Length' \, c_1 ; \break \'Length' \, a_2 ; \'Length' \, c_2 \
\"

\" \t//008. x \preceq \ep \imp x = \ep \"

\sam008.
\"\p/008.
/\'H' : x \
/264\fw ; x ; \ep \
/57 ; \'Length' \, x \
/001 ; x \
/259 \
\"

\" \t//007. a_1 \oplus c_1 = a_2 \oplus c_2 \and a_1 \preceq a_2 \and c_2 \ne \ep \imp c_1 \ne \ep \"

\sam007.
\"\p/007.
/\'H' : a_1 : c_1 : a_2 : c_2 \
/006 ; a_1 ; c_1 ; a_2 ; c_2 \
/008 ; c_2 \
\"

\" \t//003b. \ep \' begins with ' a_1 \and \ep \' begins with ' a_2 \and a_1 \preceq a_2 \imp
a_2 \' begins with ' a_1 \"

\sam003b.
\"\p/003b.
/\'H' : a_1 : a_2 \
/002 ; a_1 \
/002 ; a_2 \
/294 ; \ep \
\"

\" \at//003i.
[ x \' begins with ' a_1 \and x \' begins with ' a_2 \and a_1 \preceq a_2 \imp a_2 \' begins with ' a_1 ] \imp
[ x \oplus \u0 \' begins with ' a_1 \and x \oplus \u0 \' begins with ' a_2 \and a_1 \preceq a_2
\imp a_2 \' begins with ' a_1 ] \"

\smallskip

{\ninepoint

{\it Proof.} Suppose not. Then

\quad .1 $x\' begins with 'a_1 \andD x \' begins with ' a_2 \andD a_1 \preceq a_2 \impP a_2 \' begins with ' a_1$

\quad .2 $x \oplus \u0 \' begins with ' a_1$

\quad .3 $x \oplus \u0 \' begins with ' a_2$

\quad .4 $a_1 \preceq a_2$

\quad .5 $\neg a_2 \' begins with ' a_1$

\noi By .3 and $/290\fw;x\oplus \u0;a_2:c_2$,

\quad .6 $x\oplus\u0=a_2\oplus c_2$

\noi By $/196;a_2$,

\quad .7 $a_2\oplus \ep=a_2$

\noi Claim: $c_2\ne\ep$. Suppose not. Then

\quad\qquad .8 $c_2=\ep$

\noi By .6 and .7,

\quad\qquad .9 $x\oplus\u0=a_2$

\noi The claim is proved by .2 and .5 and .9, and .8--.9 will not be used again.

\quad .10 $c_2\ne\ep$

\noi By $/254;a_2;c_2$ and .10,

\quad .11 $\'Chop' (a_2\oplus c_2)=a_2\oplus \'Chop' \, c_2$

\noi By /243;x,

\quad .12 $\'Chop' (x\oplus \u0)=x$

\noi By .6 and .12 and .11,

\quad .13 $x=a_2\oplus \'Chop' \, c_2$

\noi By $/291;a_2;\'Chop' \,c_2$ and .13,

\quad .14 $x\' begins with ' a_2$

\noi By .2 and $/290\fw; x\oplus \u0;a_1:c_1$,

\quad .15 $x\oplus\u0=a_1\oplus c_1$

\noi By $/196;a_1$,

\quad .16 $a_1 \oplus \ep = \ep$

\noi By $/007;a_1;c_1;a_2;c_2$ and .6 and .15 and .4 and .10,

\quad .17 $c_1\ne\ep$

\noi By $/254;a_1;c_1$ and .17,

\quad .18 $\'Chop' (a_1\oplus c_1)=a_1\oplus \'Chop' \, c_1$

\noi By $/291; a_1;\'Chop' c_1$ an .18 and .12 and .15,

\quad .19 $x \' begins with ' a_1$

\noi QEA by .1 and .19 and .14 and .4 and .5.}

\" \at//003j.
[ x \' begins with ' a_1 \and x \' begins with ' a_2 \and a_1 \preceq a_2 \imp a_2 \' begins with ' a_1 ] \imp
[ x \oplus \u1 \' begins with ' a_1 \and x \oplus \u1 \' begins with ' a_2 \and a_1 \preceq a_2
\imp a_2 \' begins with ' a_1 ] \"

\smallskip

The proof is entirely similar: replace each $\u0$ by $\u1$ and /243 by /244.

\" \at/003. x \' begins with ' \, a_1 \and x \' begins with ' a_2 \and a_1 \preceq a_2 \imp
a_2 \' begins with ' a_1 \"

\" \t//010. x \' begins with ' y \and y \' begins with ' x \imp x = y \"

\sam010.
\"\p/010.
/\'H' : x : y \
/290\fw ; x ; y : c \
/290\fw ; y ; x : d \
/183 ; x ; d ; c \
/196 ; x \
/199 ; x ; \ep ; d \oplus c \
/204 ; d ;c \
\"

\" \t//011. x \' ends with ' y \and y \' ends with ' x \imp x = y \"

\sam011.
\"\p/011.
/\'H' : x : y \
/288\fw ; x ; y : c \
/288\fw ; y ; x : d \
/194 ; x \
/183 ; c ; d ; x \
/197 ; \ep ; x ; c \oplus d \
/204 ; c ; d \
\"

\" \d//012. \'last-bit' ( x , b ) \iff [ x = \'Chop' x \oplus \u0 \and b = \u0 ] \or
[ x = \'Chop' x \oplus \u1 \and b = \u1 ] \or [ x =\ep \and b = \ep ] \"

\" \t//013. \'last-bit' ( x , b_1 ) \and \'last-bit' ( x , b_2 ) \imp b_1 = b_2 \"

\sam013.
\"\p/013.
/\'H' :x : b_1 : b_2 \
/012\fw ; x ; b_1 \
/012\fw ; x ; b_2 \
/253 ; x \
/? x \ne \ep \
\"

\ifx\wrapperUsed\nil
\def\thisEnd{\end}
\else
\def\thisEnd{\relax}
\fi

\thisEnd

\long\def\eat#1{\relax}
\eat{    REST IS IGNORED!!!
\"stop\"

{\it \#15. Strings of length two are called marks. Every string of even length is uniquely a concatenation
of marks. Introduce the notion of even occurrence expressing that one concatenation of marks occurs in
another. A list is a concatenation of marks ending with the mark\/ $\u0 \oplus \u1 $. An atomic list is
a minimal list, and every list is uniquely a concatenation of atomic lists.}

\" \d//279. x \' has even length' \iff \'Parity' \, \'Length' \, x = 0 \hfill\star \"

\" \t//280. x \' has even length' \and y \' has even length' \imp x \oplus y \' has even length' \"

\sam280.
\"\p/280.
/\'H' : x : y \
/279\fw ; x \
/279\fw ; y \
/260 ; x ; y \
/210 ; \'Length' \, x ; \'Length' \, y \
/279\bw ; x \oplus y \

\"

\" \t//281. x \' has even length' \and x \oplus y \' has even length' \imp y \' has even length' \"

\sam281.
\"\p/281.
/\'H' : x : y \
/279\fw ; x \
/279\fw ; x \oplus y \
/260 ; x ; y \
/210 ; \'Length' \, x ; \'Length' \, y \
/279\bw ; y \
\"

\" \t//282. x \oplus  y \' has even length' \and y \' has even length' \imp x \' has even length' \"

\sam282.
\"\p/282.
/\'H' : x : y \
/279\fw ; x \oplus y \
/279\fw ; y \
/260 ; x ; y \
/98 ; \'Length' \, x ; \'Length' \, y \
/210 ; \'Length' \, y ; \'Length' \, x \
/279\bw ; x \
\"

\" \t//283. \u0 \oplus \u0 \' has even length' \and
\u0 \oplus \u1 \' has even length' \and
\u1 \oplus \u0 \' has even length' \and
\u1 \oplus \u1 \' has even length' \"

\sam283.
\"\p/283.
/\'H' \
/278 \
/260 ; \u0 ; \u0 \
/260 ; \u0 ; \u1 \
/260 ; \u1 ; \u0 \
/260 ; \u1 ; \u1 \
/279\bw ; \u0 \oplus  \u0 \
/279\bw ; \u0 \oplus  \u1 \
/279\bw ; \u1 \oplus  \u0 \
/279\bw ; \u1 \oplus  \u1 \
/208 \
/115 \
/116 \
/12 ; \'S' 0 ; 0 \
/12 ; \'S' \'S' 0 \
\"

\" \d//284. x \' occurs evenly in ' y \iff \exists a \preceqQ y \exists c \preceqQ y
[ a \' has even length' \and x \' has even length' \and c \' has even length' \and
y = a \oplus x \oplus c ] \hfill\star \"

\smallskip

Call $\u0\oplus\u0$, $\u1\oplus\u1$, $\u0\oplus\u1$, and $\u1\oplus\u0$ \ii marks.. They are the
strings of length two. Then \hbox{``$x \' has even length'$''} expresses that $x$ is a concatenation of marks,
and ``$x \' occurs evenly in~'y$'' expresses that $x$~occurs as a concatenation of marks in~$y$.
The alphabet of two bits has been expanded to the alphabet of four marks.

\" \t//285. a \' has even length' \and x \' has even length' \and c \' has even length' \imp
x \' occurs evenly in ' a \oplus x \oplus c \"

\sam285.
\"\p/285.
/\'H' : a : x : c \
/284\bw ; x ; a \oplus x \oplus c ; a ; c \
/277 ; a ; x \oplus  c \
/183 ; a ; x ; c \
/277 ; a \oplus x ; c \
/5 ; a \oplus x \oplus c \
\"

\" \t//286. x \' occurs evenly in ' y \and z \' has even length' \imp x \' occurs evenly in ' y \oplus z \"

\sam286.
\"\p/286.
/\'H' : x : y : z \
/284\fw ; x ; y : a : c \
/183 ; a \oplus x ; c ; z \
/183 ; a ; x ; c \oplus z \
/183 ; a ; x ; c \
/280 ; c ; z \
/285 ; a ; x ; c \oplus z \
\"

\" \t//287. x \' occurs evenly in ' y \and z \' has even length' \imp x \' occurs evenly in ' z \oplus y \"

\sam287.
\"\p/287.
/\'H' : x : y : z \
/284\fw ; x ; y : a : c \
/183 ; z ; a ; x \oplus c \
/280 ; z ; a \
/285 ; z \oplus a ; x ; c \
\"

\" \d//295. L \' is a list' \iff L \' has even length' \and L \' ends with ' \u0 \oplus \u1 \hfill\star \"

\" \t//296. L_1 \' is a list' \and L_2 \' is a list' \imp L_1 \oplus L_2 \' is a list' \"

\sam296.
\"\p/296.
/\'H' : L_1 : L_2 \
/295\fw ; L_1 \
/295\fw ; L_2 \
/280 ; L_1 ; L_2 \
/293 ; L_2 ; \u0 \oplus \u1 ; L_1 \
/295\bw ; L_1 \oplus L_2 \
\"


\" \d//297. L_1 \' is a sublist of ' L_2 \iff
L_1 \' is a list' \and L_2 \' is a list' \and
[ \u0 \oplus \u1 \oplus L_1 \' occurs evenly in ' L_2 \or L_2 \' begins with ' L_1 ] \hfill\star \"

\" \t//298. L_1 \' is a sublist of ' L_2 \and L_3 \' is a list' \imp L_1 \' is a sublist of ' L_2 \oplus L_3 \"

\smallskip

If $L_2 \' begins with ' L_1$, so does $L_2\oplus L_3$. If $\u0 \oplus \u1 \oplus L_1 \' occurs evenly in ' L_2$,
then it does in $L_2\oplus L_3$.

\sam298.
\"\p/298.
/\'H' : L_1 : L_2 : L_3 \
/297\fw ; L_1 ; L_2 \
/295\fw ; L_3 \
/296 ; L_2 ; L_3 \
/297\bw ; L_1 ; L_2 \oplus L_3 \
/292 ; L_2 ; L_1 ; L_3 \
/286 ; \u0 \oplus \u1 \oplus L_1 ; L_2 ; L_3 \
\"

\" \t//299. L_1 \' is a sublist of ' L_2 \and L_3 \' is a list' \imp L_1 \' is a sublist of ' L_3 \oplus L_2 \"

\smallskip

If $\u0 \oplus \u1 \oplus L_1 \' occurs evenly in ' L_2$, then is does in $L_3\oplus L_2$. If
$L_2 \' begins with ' L_1$, then there is a~$c$ such that $L_2=L_1\oplus c$. Since
$L_3 \' ends with ' \u0 \oplus \u1 $, there is an~$a$ such that $L_3=a\oplus \u0 \oplus \u1 $.
Verify that $a$ and~$c$ have even length, and re-associate to find $L_3\oplus L_2=
a\oplus (\u0 \oplus \u1 \oplus L_1)\oplus c$.

\sam299.
\"\p/299.
/\'H' : L_1 : L_2 : L_3 \
/297\fw ; L_1 ; L_2 \
/295\fw ; L_1 \
/295\fw ; L_2 \
/295\fw ; L_3 \
/297\bw ; L_1 ; L_3 \oplus L_2 \
/296 ; L_3 ; L_2 \
/287 ; \u0 \oplus \u1 \oplus L_1 ; L_2 ; L_3 \
/288\fw ; L_3 ; \u0 \oplus \u1 : a \
/283 \
/282 ; a ; \u0 \oplus \u1 \
/290\fw ; L_2 ; L_1 : c \
/281 ; L_1 ; c \
/280 ; \u0 \oplus \u1 ; L_1 \
/183 ; \u0 ; \u1 ; L_1 \
/183 ; a ; \u0 \oplus \u1 ; L_1 \oplus c \
/183 ; \u0 \oplus \u1 ; L_1 ; c \
/285 ; a ; \u0 \oplus \u1 \oplus L_1 ; c \
\"

\" \d//300. \'even-occurrence' ( a , x , c , y ) \iff a \' has even length' \and x \' has even length' \and
c \' has even length' \and y = a \oplus x \oplus c \"

}


%% file: NelsonarXivBalrogMacros.tex
\mathcode`: = 58
\mathcode`; = 59

\font\tenmsa=msam10

\font\ninerm=cmr9
\font\ninebf=cmbx9
\font\nineit=cmti9
\font\ninei=cmmi9
\font\ninett=cmtt9
\font\eighti=cmmi8
\font\eightrm=cmr8
\font\sevenrm=cmr7
\font\eightbf=cmbx8
\font\eightit=cmti8

\font\eighttt=cmtt8
\font\eightsy=cmsy8

\font\ninesy=cmsy9
\font\eightsy=cmsy8
\font\nineit=cmti9
\font\eightit=cmti8
\font\sixrm=cmr6
\font\sixbf=cmbx6
\font\sixi=cmmi6
\font\sixsy=cmsy6

\font\eightsans=cmss8

\def\ninepoint{\def\rm{\fam0\ninerm}%
   \textfont0=\ninerm \scriptfont0=\sixrm \scriptscriptfont0=\fiverm
   \textfont1=\ninei \scriptfont1=\sixi \scriptscriptfont1=\fivei
   \textfont2=\ninesy \scriptfont2=\sixsy \scriptscriptfont2=\fivesy
   \textfont3=\tenex \scriptfont3=\tenex \scriptscriptfont3=\tenex
   \textfont\itfam=\nineit \def\it{\fam\itfam\nineit}%
   \textfont\bffam=\ninebf \scriptfont\bffam=\sixbf
    \scriptscriptfont\bffam=\fivebf \def\bf{\fam\bffam\ninebf}%
   \normalbaselineskip=11pt
   \setbox\strutbox=\hbox{\vrule height8pt depth 3pt width0pt}%
   \let\big=\ninebig \normalbaselines\rm}
\def\eightpoint{\def\rm{\fam0\eightrm}%
   \textfont0=\eightrm \scriptfont0=\sixrm \scriptscriptfont0=\fiverm
   \textfont1=\eighti \scriptfont1=\sixi \scriptscriptfont1=\fivei
   \textfont2=\eightsy \scriptfont2=\sixsy \scriptscriptfont2=\fivesy
   \textfont3=\tenex \scriptfont3=\tenex \scriptscriptfont3=\tenex
   \textfont\itfam=\eightit \def\it{\fam\itfam\eightit}%
   \textfont\bffam=\eightbf \scriptfont\bffam=\sixbf
    \scriptscriptfont\bffam=\fivebf \def\bf{\fam\bffam\eightbf}%
   \normalbaselineskip=11pt
   \setbox\strutbox=\hbox{\vrule height8pt depth 3pt width0pt}%
   \let\big=\eightbig \normalbaselines\rm}
\catcode`@=11
\def\ninebig#1{{\hbox{$\textfont0=\tenrm\textfont2=\tensy
  \left#1\vbox to7.25pt{}\right.\n@space$}}}
\catcode`@=12



\def\yy{\catcode`_=12\catcode`/=12\catcode`~=12\catcode`_=12}
\def\zz{\catcode`/=\active\catcode`~=13\catcode`_=8}
\yy
\def\sam#1.{\smallskip\par{\quad}{{\yy\eightit\pdfklink{Proof.}{http://math.princeton.edu/~nelson/proof/#1.pdf}}}\zz\hskip8pt }
\zz
\yy
\def\rf#1;{\pdfklink{#1;}{http://math.princeton.edu/~nelson/ref/#1.pdf}\zz}
\zz
\yy
\def\Rf#1{\pdfklink{#1}{http://math.princeton.edu/~nelson/ref/#1.pdf}\zz}
\zz


\def\red#1{\rgbo{1 0 0}{#1}}

\yy\zz


\let\doublequote="
\def\"{}
\def\slash{\char'57}
\catcode`/=13
\def/{}
\def\a#1.#2\par{\smallbreak\par\noindent a$#1$.\quad$#2$}
\def\ar#1.#2\par{\smallbreak\par\noindent r$#1$.\quad$#2$}
\def\at#1.#2\par{\smallbreak\par\noindent t$#1$.\quad$#2$}
\def\at#1.#2\par{\smallbreak\par\noindent t$#1$.\quad$#2$}
\def\t#1.#2\par{\smallbreak\par\noindent t$#1$.\quad$#2$}
\def\tu#1.#2\par{\medbreak\par\noindent t$#1$.~({\eightrm UC})\quad$#2$}
\def\te#1.#2\par{\medbreak\par\noindent t$#1$.~({\eightrm EC})\quad$#2$}
\def\d#1.#2\par{\smallbreak\par\noindent d$#1$.\quad$#2$}
\def\de#1.#2\par{\smallbreak\par\noindent e$#1$.\quad$#2$}
\def\p/#1.#2\"{\begingroup\baselineskip=9pt\everymath={\scriptstyle}\def\rm{\sevenrm}\def\ {\hskip 10pt plus 2pt minus 10pt}$\!\!#2$\par\endgroup}
\newif\ifproof

\let\acute=\'
\def\'#1'{{}$\rm#1${}}
\def\iff{\allowbreak\hskip8pt plus2pt minus3pt\mathrel{\hbox{$\leftrightarrow$}}\hskip8pt plus2pt minus3pt}
\def\iffF{\allowbreak\leftrightarrow}
\def\diff{\allowbreak\mathrel{\dot\leftrightarrow}}
\def\imp{\allowbreak\hskip8pt plus2pt minus3pt\mathrel{\hbox{$\rightarrow$}}\hskip8pt plus2pt minus3pt}
\def\impP{\allowbreak\rightarrow}
\def\dimp{\allowbreak\mathrel{\dot\rightarrow}}

\def\andD{\allowbreak\mathrel{\hskip1pt\&\hskip1pt}}
\def\dand{\allowbreak\mathrel{\hskip-1pt{\dot{\hskip3pt\&}}\hskip1pt}}

\def\and{\allowbreak\hskip7pt plus2pt minus3pt\mathrel{\&}\hskip7pt plus2pt minus3pt}
\def\orR{\allowbreak\kern1pt\lor\kern1.25pt}
\def\dor{\allowbreak\mathrel{\kern1pt\dot\lor\kern1.25pt}}
\def\dotor{\allowbreak\mathrel{\kern1pt\dot\lor\kern1.25pt}}
\def\or{\allowbreak\hskip7pt plus2pt minus3pt\mathrel{\hbox{$\vee$}}\hskip7pt plus2pt minus3pt}
\def\deq{\mathrel{\dot=}}

\def\dneg{\dot\Neg}

\def\leE{{\scriptstyle\le}}
\def\lt{<}
\def\ltT{{\scriptstyle<}}
\def\preceqQ{{\scriptstyle\preceq}}

\newcount\firstpage
\newcount\footno
\newcount\secno
\newtoks\chaptername
\newtoks\sectionname
\newwrite\tabcon
\newif\ifindexmode

\predisplaypenalty=0

\let\ep=\epsilon
\let\epsilon=\varepsilon

\def\cite#1{[{\ninerm Ne~Ch.~#1}]}
\def\ii#1.{{\it#1}}

\let\phi=\varphi

\let\hyphen=\-
\def\-{\relax\hyphen}
\def\ro#1{{\rm #1}}

\def\bul{\hfill\hbox{\tenmsa\char3}}

\smallskip
\def\Pr{\smallskip\par P{\eightrm ROPOSITION}}
\def\med#1.#2\par{\smallskip\par\noindent #1.\quad$#2$\par}

\def\pf{\par\smallskip\leavevmode\hbox{\it Proof. }}

\let\Neg=\neg
\def\neg{\Neg\;}
\def\noin{\smallskip\par\noindent}
\def\noi{\noindent}
\def\nn#1.{$(#1)$}

\def\nok// #1.{\par\smallskip\par\noindent(#1)\quad}
\def\Nok/ #1.{\par\smallskip\par\noindent(#1)\quad}

\def\pphi|#1|{\phi(#1)}
\def\pphiz|#1|{\phi^0(#1)}
\def\pphio|#1|{\phi^1(#1)}
\def\pphit|#1|{\phi^2(#1)}
\def\ppphiz|#1|{\phi^{(0)}(#1)}
\def\ppphio|#1|{\phi^{(1)}(#1)}
\def\ppphit|#1|{\phi^{(2)}(#1)}
\long\def\omit#1{}


\let\Norm=\|

\def\|{{\scriptstyle\le}}

\def\0{\hbox{\rm o}}
\def\1{\hbox{\i}}
\def\2{\hbox{\eighttt|}}
\def\3{\hbox{\eighttt @}}

\let\lbracket=[
\let\rbracket=]
\def\[{\,\big\lbracket\,}
\let\lbrace=\{
\let\rbrace=\}


\def\Qu[#1]{\langle#1\rangle}


\def\bw{{^{\scriptscriptstyle<\!\!-}}}
\def\fw{{^{\scriptscriptstyle-\!\!>}}}
\def\tauu_#1{\tau_{\scriptscriptstyle#1}}

\let\lblank=\ %

\def\true{\hbox{\eightsans T}}
\def\false{\hbox{\eightsans F}}
\def\minus{-}
\def\down#1{^{\phantom'}_{#1}}
\def\Pair (#1,#2){\langle#1,#2\rangle}

\def\up{\mathbin{\uparrow}}

\def\u#1{\underline#1}

%% file: BussTao_Afterword_arXiv.tex
\ifx\wrapperUsed\nil
\font\sambigrm=cmr10 scaled \magstep2
\fi
\font\ninett=cmtt9
\font\ninerm=cmr9
\font\ninebf=cmbx9
\font\nineit=cmti9
\font\ninei=cmmi9
\font\ninett=cmtt9
\font\eighti=cmmi8
\font\eightrm=cmr8
\font\eightbf=cmbx8
\font\eightit=cmti8

\font\eighttt=cmtt8
\font\eightsy=cmsy8

\font\ninesy=cmsy9
\font\eightsy=cmsy8
\font\nineit=cmti9
\font\eightit=cmti8
\font\sixrm=cmr6
\font\sixbf=cmbx6
\font\sixi=cmmi6
\font\sixsy=cmsy6
\def\ninepoint{\def\rm{\fam0\ninerm}%
   \textfont0=\ninerm \scriptfont0=\sixrm \scriptscriptfont0=\fiverm
   \textfont1=\ninei \scriptfont1=\sixi \scriptscriptfont1=\fivei
   \textfont2=\ninesy \scriptfont2=\sixsy \scriptscriptfont2=\fivesy
   \textfont3=\tenex \scriptfont3=\tenex \scriptscriptfont3=\tenex
   \textfont\itfam=\nineit \def\it{\fam\itfam\nineit}%
   \textfont\bffam=\ninebf \scriptfont\bffam=\sixbf
    \scriptscriptfont\bffam=\fivebf \def\bf{\fam\bffam\ninebf}%
   \normalbaselineskip=11pt
   \setbox\strutbox=\hbox{\vrule height8pt depth 3pt width0pt}%
   \let\big=\ninebig \normalbaselines\rm}

\ifx\wrapperUsed\relax
\magnification=\magstep1
\fi
\hsize=5.0in
\hoffset=0.25in

\def\bibindent{\vskip0.5em\par\noindent\hangindent=1em\hangafter=1}
\def\biblabel#1{#1\hskip4pt\ignorespaces}
\def\BussNelsonPred{[1]}
\def\FarisNelson{[2]}
\def\KritchmanRazsurprise{[3]}
\def\NelsonIST{[4]}
\def\NelsonPredArith{[5]}
\def\YesseninVolpin{[6]}

\def\IDeltazero{{\rm I} \Delta_0}

\
\vskip1cm

\centerline{
\sambigrm  Afterword on two works by Ed Nelson}

\vskip1em

\centerline{
Sam Buss
~~ and ~~ Terence Tao
}

\vskip 2em
Two of Ed Nelson's unfinished papers, {\it Elements} (dated March 12, 2013)
and {\it Inconsistency of Primitive Recursive Arithmetic} (undated, also known
as the {\it Balrog} paper)
have the goal of proving the inconsistency of
number theory.
As Nelson writes in {\it Elements}, ``The aim of this work is to show that
contemporary mathematics, including Peano arithmetic, is inconsistent ...''.

Neither of these papers have been circulated before in their
current forms.  An earlier version of {\it Elements} was circulated
in 2011, but was found to have problems in its treatment of
proofs generated by Chaitin machines.  The new 2013 version
uses a similar approach, but gives a much more detailed
explanation of the planned proof, and it
handles Chaitin machines differently so as to address
the earlier problems.

Nelson's remarkable program to establish the inconsistency
of Peano arithmetic was intertwined with his development
of Internal Set Theory~\NelsonIST{}
and especially Predicative Arithmetic~\NelsonPredArith{}.
Predicative Arithmetic is a constructive fragment
of arithmetic, and Nelson's development of Predicative
Arithmetic was inspired in part
by Yessenin-Volpin's ultra-intuitionistic
set theory~\YesseninVolpin{}.
The mathematical content of Predicate Arithmetic
is closely tied to theories of bounded
arithmetic such as $\IDeltazero$, $\IDeltazero+\Omega_1$,
${\rm S} ^i_2$ and~${\rm T}^i_2$.
Indeed, Nelson~\NelsonPredArith{}
independently discovered some of the important
tools for bounded arithmetic, including the technique of
speeding up induction on cuts and the local interpretability
of predicative arithmetic and bounded arithmetic in
Robinson's theory~$Q$.  Nelson's Predicative Arithmetic
was also influential for the definition by one of us (Buss)
of the theories ${\rm S}^i_2$ and~${\rm T}^i_2$ of
bounded arithmetic, including notably the use of
Nelson's smash function.

The {\it Elements} manuscript gives a detailed,
high-level outline of Nelson's plan for a proof
of the inconsistency of Peano arithmetic (and
primitive recursive arithmetic).  One
of the principal tools is a novel use
of a recent proof
by Kritchman and Raz~\KritchmanRazsurprise{}
of G\"odel's second incompleteness theorem
based on the ``surprise examination''.
Nelson also uses Kolmogorov complexity and
techniques from cut-elimination.  The detailed plan
of the inconsistency proof
is outlined as Steps 1-17 in section~9 near the
end of {\it Elements}.  The plan first discusses a
system~$S^*$ which is a predicative theory including
bounded induction and which is strong enough to express
concepts about metamathetical concepts,
Chaitin machine computation, and Kolmogorov complexity.
Step~7 introduces a finitary theory~$\cal F$; the details
of the system~$\cal F$ are not fully specified, but it needs
to be able to formalize cut-elimination or normalization.
Thus it seems that $\cal F$ can be taken to be, for instance,
$\IDeltazero+{\rm superexp}$ or ${\rm I} \Sigma_1$.
The heart of the argument is reached in Step~16.
Unfortunately, the argument becomes very uncertain here.
Nelson argues that $S^*$ disproves a sequence of statements:
first $A_{\kappa, 0}$, then $A_{\kappa,1}$, etc.,
up through $A_{\kappa,I}$.
The base case that $S^*$ disproves $A_{\kappa,0}$ is fine,
but the later stages are unclear.
It seems that the disproof of $A_{\kappa,\delta}$
requires an assumption that $S^*$ is consistent.
The reason for this is that the Chaitin machine cannot be
given the value of~$\delta$ as an input since the
Kolmogorov complexity of~$\delta$ may not be sufficiently
below that of~$\kappa$. The only alternative to
explicitly specifying~$\delta$ that we can think of,
is for the Chaitin machine to first search
for the $S^*$~proof that ``$\lnot A_{\kappa,\delta+1}$''
and then also wait until $\delta$ many strings are
found to have Chaitin complexity less than~$\kappa$.
This however assumes that $S^*$ is consistent.
Of course, $S^*$ does not prove its own consistency.
Perhaps Nelson had a different argument in mind,
but this is our best attempt to flesh out his arguments.
At any rate,
Nelson was apparently aware of the potential problem here,
since he earlier discusses the need for a system to prove
the ``consistency of its own arithmetization''.

In the spirit of a quote by Carl Sagan, ``Extraordinary claims
require extraordinary evidence'', Nelson planned to fulfill his
inconsistency proof by exhibiting a fully formal,
computer-verified derivation of a contradiction.  That is,
he planned not to prove that there is a proof of contradiction,
but to actually exhibit an explicit proof of a contradiction.  The
first steps of this are carried out at the end of {\it Elements},
and it is even further pursued in {\it Balrog}.
The {\it Balrog} manuscript is still incomplete, as only
six sections are complete, and at least ten sections were planned.
The {\it Balrog} manuscript is in essence
a formalization of the ``bootstrapping''
of predicative arithmetic in the spirit
of~\NelsonPredArith{}.  A remarkable
feature of {\it Balrog} is that
proofs of theorems are indicated in a terse fashion that
permits a Perl program, called {\it qea}, to automatically verify the
proofs.  For instance, Theorems {\it 13b.}\
and {\it 13i.} of~{\it Balrog},
and their proofs, are typeset with
the TeX code
\vskip1em

{\ninett
\noindent
\string\" \string\t\string/\string/13b. 0 + 0 = 0 + 0 \string\"
\vskip1em

\noindent
\string\sam13b.

\noindent
\string\"\string\p\string/13b.

\noindent
\string/\string\'H' \string\ {}

\noindent
\string/5 ; 0 + 0 \string\ {}

\noindent
\string\"
\vskip1em

\noindent
\string\" \string\t\string/13i. x + 0 = 0 + x \string\imp{} \string\'S' x + 0 = 0 + \string\'S' x \string\"
\vskip1em

\noindent
\string\sam13i.

\noindent
\string\"\string\p\string/13i.

\noindent
\string/\string\'H' : x \string\ {}

\noindent
\string/12 ; \string\'S' x \string\ {}

\noindent
\string/12 ; 0 ; x \string\ {}

\noindent
\string/12 ; x \string\ {}

\noindent
\string\"

}
\vskip 1em

These indicate that {\it 13b.}\ is proved by
substituting 0+0 for~$x$ in axiom {\it a5.}, and that
{\it 13i.}\ is proved by using
definition {\it r12.}\ three times,
first substituting
$Sx$ for~$x$, then $0$ and~$x$
for $x$ and~$y$, and finally
$x$ for~$x$.  After these substitutions, the
desired conclusions follow propositionally 
from equality axioms.  The
{\it qea} system then automatically generated
an expanded proof; the expanded proof was produced
as a TeX file, and automatically converted to PDF.
An example is shown in an appendix to the
{\it Balrog} manuscript posted to the {\it arXiv}.

We of course believe that Peano
arithmetic is consistent; thus we do not expect
that Nelson's project can be completed according to
his plans.  Nonetheless, there is much new
in his papers that is of potential mathematical,
philosophical and computational interest.  For this reason,
they are being posted to the {\it arXiv}.  Two aspects
of these papers seem particularly useful.  The first aspect is
the novel use of the ``surprise examination'' and Kolmogorov
complexity; there is some possibility that similar techniques might
lead to new separation results for fragments of
arithmetic.  The second aspect is Nelson's automatic proof-checking
via TeX and {\it qea}.  This is highly
interesting and provides
a novel method of integrating human-readable
proofs with computer verification of proofs.
\medskip

The reader interested in further discussion of
Nelson's Predicative Arithmetic can
consult the mostly-survey article~\BussNelsonPred{}.
The volume~\FarisNelson{} contains papers about
many other aspects of Nelson's wide-ranging research.
Other works by Nelson are available at
{\tt math.princeton.edu\slash$\sim$nelson}, including a number
of philosophical works.

\vskip0.5em

{\ninepoint
\bibindent
\biblabel{\BussNelsonPred}
{S.~R.\ Buss}, {\it Nelson's work on logic and foundations and other
  reflections on foundations of mathematics}, in Diffusion, Quantum Theory, and
  Radically Elementary Mathematics, Princeton University Press, 2006,
  pp.~183--208. Edited by W.\ Faris.
\bibindent
\biblabel{\FarisNelson}
{W.~G.\ Faris}, ed., {\it Diffusion, Quantum Theory, and Radically
  Elementary Mathematics}, Mathematical Notes, \#47, Princeton University
  Press, 2006.
\bibindent
\biblabel{\KritchmanRazsurprise}
{S.~Kritchman and R.~Raz}, {\it The surprise examination and the second
  incompleteness theorem}, Notices of the American Mathematical Society, 57
  (2010), pp.~1454--1458.
\bibindent
\biblabel{\NelsonIST}
{E.~Nelson}, {\it Internal set theory: A new approach to nonstandard
  analysis}, Bulletin of the American Mathematical Society, 83 (1977),
  pp.~1165--1198.
\bibindent
\biblabel{\NelsonPredArith}
{E.~Nelson}, {\it {P}redicative {A}rithmetic}, Princeton University Press,
  1986.
\bibindent
\biblabel{\YesseninVolpin}
{A.~S.\ Yessenin-Volpin}, {\it The ultra-intuitionistic criticism and the
  antitraditional program for foundations of mathematics}, in Intuitionism and
  Proof Theory, A.~Kino, J.~Myhill, and R.~E. Vesley, eds., North-Holland,
  1970, pp.~1--45.
}

\bigskip
\hfill
{Sam Buss}
\vskip 0.5ex

\hfill
{Terence Tao}
\vskip 2ex

\hfill {\it
August 31, 2015}

\ifx\wrapperUsed\nil
\def\thisEnd{\end}
\else
\def\thisEnd{\relax}
\fi

\thisEnd

%% file: Appendix13i.tex
\vskip1cm

\centerline{
\sambigrm  Appendix}

\vskip1em

\centerline{
Example: Expanded proof of Balrog 13i}

\vskip 2em
Nelson's {\it qea} proof system consists
of a Perl script
which reads the TeX source code, and checks the proof correctness
and generates an expanded version of the proof.
The {\it Balrog} and {\it Elements} documents
contained active hyperlinks, in a blue font,
to these expanded proofs.

As an example, the expanded version
of the proof of {\it 13i.}\ in {\it Balrog} as
generated by {\it qea} is shown below.

\vskip 4em

\hrule

\vskip4em

\def\red#1{\hbox{\textRed{$#1$}\textBlack}}

\parindent=0pt

\raggedright

	\overfullrule=0pt \parskip=\medskipamount
	\centerline{\bf Proof of Theorem 13i}

\bigskip

	\quad The theorem to be proved is

 $  x + 0 = 0 + x \imp \'S' x + 0 = 0 + \'S' x  $

	\quad Suppose the theorem does not hold. Then, with the variables held fixed,

 $(\ro H)\quad [ [ ( x + 0 ) = ( 0 + x ) ] \and [ \neg ( ( \'S' x ) + 0 ) = ( 0 + ( \'S' x ) ) ] ]  $

	\bigskip

 \centerline{\bf Special cases of the hypothesis and previous results:}\bigskip

\leavevmode\phantom00:\quad$0 + x = x + 0 $\qquad from\quad $ /\ro H : x $

\leavevmode\phantom01:\quad$\neg ( \'S' x ) + 0 = 0 + ( \'S' x ) $\qquad from\quad $ /\ro H : x $

\leavevmode\phantom02:\quad$( \'S' x ) + 0 = \'S' x  $\qquad from\quad $ \yy\pdfKlink{12}{http://math.princeton.edu/~nelson/ref/12.pdf}\zz  ; \'S' x  $

\leavevmode\phantom03:\quad$\'S' ( 0 + x ) = 0 + ( \'S' x ) $\qquad from\quad $ \yy\pdfKlink{12}{http://math.princeton.edu/~nelson/ref/12.pdf}\zz  ; 0 ; x $

\leavevmode\phantom04:\quad$x + 0 = x $\qquad from\quad $ \yy\pdfKlink{12}{http://math.princeton.edu/~nelson/ref/12.pdf}\zz  ; x $

\bigskip

\centerline{\bf Equality substitutions:}

\bigskip

\leavevmode\phantom05:\quad$\neg 0 + x = x + 0 \or \neg \'S' ( \red{0 + x} ) = 0 + ( \'S' x ) \or \'S' ( \red{x + 0} ) = 0 + ( \'S' x ) $\medskip

\leavevmode\phantom06:\quad$\neg ( \'S' x ) + 0 = \'S' x  \or \red{( \'S' x ) + 0} = 0 + ( \'S' x ) \or\neg \red{\'S' x} = 0 + ( \'S' x ) $\medskip

\leavevmode\phantom07:\quad$\neg x + 0 = x \or \neg \'S' ( \red{x + 0} ) = 0 + ( \'S' x ) \or \'S' ( \red{x} ) = 0 + ( \'S' x ) $\medskip

\bigskip

\centerline{\bf Inferences:}\openup1\jot

\bigskip

\leavevmode\phantom08:\quad$\neg \'S' ( 0 + x ) = 0 + ( \'S' x ) \or \'S' ( x + 0 ) = 0 + ( \'S' x ) $\qquad by \hfill\break\phantom.\qquad0: $\red{0 + x = x + 0 }$\quad \hfill\break\phantom.\qquad5: $\red{\neg 0 + x = x + 0 }\or \neg \'S' ( 0 + x ) = 0 + ( \'S' x ) \or \'S' ( x + 0 ) = 0 + ( \'S' x ) $

\leavevmode\phantom09:\quad$\neg ( \'S' x ) + 0 = \'S' x  \or \neg 0 + ( \'S' x ) = \'S' x  $\qquad by \hfill\break\phantom.\qquad1: $\red{\neg ( \'S' x ) + 0 = 0 + ( \'S' x ) }$\quad \hfill\break\phantom.\qquad6: $\neg ( \'S' x ) + 0 = \'S' x  \or \red{( \'S' x ) + 0 = 0 + ( \'S' x ) }\or \neg 0 + ( \'S' x ) = \'S' x  $

10:\quad$\neg 0 + ( \'S' x ) = \'S' x  $\qquad by \hfill\break\phantom.\qquad2: $\red{( \'S' x ) + 0 = \'S' x  }$\quad \hfill\break\phantom.\qquad9: $\red{\neg ( \'S' x ) + 0 = \'S' x  }\or \neg 0 + ( \'S' x ) = \'S' x  $

11:\quad$\'S' ( x + 0 ) = 0 + ( \'S' x ) $\qquad by \hfill\break\phantom.\qquad3: $\red{\'S' ( 0 + x ) = 0 + ( \'S' x ) }$\quad \hfill\break\phantom.\qquad8: $\red{\neg \'S' ( 0 + x ) = 0 + ( \'S' x ) }\or \'S' ( x + 0 ) = 0 + ( \'S' x ) $

12:\quad$\neg \'S' ( x + 0 ) = 0 + ( \'S' x ) \or 0 + ( \'S' x ) = \'S' x  $\qquad by \hfill\break\phantom.\qquad4: $\red{x + 0 = x }$\quad \hfill\break\phantom.\qquad7: $\red{\neg x + 0 = x }\or \neg \'S' ( x + 0 ) = 0 + ( \'S' x ) \or 0 + ( \'S' x ) = \'S' x  $

13:\quad$\neg \'S' ( x + 0 ) = 0 + ( \'S' x ) $\qquad by \hfill\break\phantom.\qquad10: $\red{\neg 0 + ( \'S' x ) = \'S' x  }$\quad \hfill\break\phantom.\qquad12: $\neg \'S' ( x + 0 ) = 0 + ( \'S' x ) \or \red{0 + ( \'S' x ) = \'S' x  }$

14:\quad$QEA$\qquad by \hfill\break\phantom.\qquad11: $\red{\'S' ( x + 0 ) = 0 + ( \'S' x ) }$\quad \hfill\break\phantom.\qquad13: $\red{\neg \'S' ( x + 0 ) = 0 + ( \'S' x ) }$

\bye